



\font\bmit=cmmib10
\font\bmits=cmmib7
\def\boldmath#1{\mathchoice{\hbox{\bmit #1}}%
{\hbox{\bmit #1}}{\hbox{\bmits #1}}{\hbox{\bmits #1}}}
\font\sc=cmcsc10
\font\LARGEbf=cmbx10 scaled \magstep4
\font\largebf=cmbx10 scaled \magstep2


\catcode`\@=11

\font\tenmsa=msam10
\font\sevenmsa=msam7
\font\fivemsa=msam5
\newfam\msafam
\textfont\msafam=\tenmsa
\scriptfont\msafam=\sevenmsa
\scriptscriptfont\msafam=\fivemsa
\font\tenmsb=msbm10
\font\sevenmsb=msbm7
\font\fivemsb=msbm5
\newfam\msbfam
\textfont\msbfam=\tenmsb
\scriptfont\msbfam=\sevenmsb
\scriptscriptfont\msbfam=\fivemsb

\def\hexnumber@#1{\ifcase#1 0\or1\or2\or3\or4\or5\or6\or7\or8\or9\or
	A\or B\or C\or D\or E\or F\fi }

\font\teneufm=eufm10
\font\seveneufm=eufm7
\font\fiveeufm=eufm5
\newfam\eufmfam
\textfont\eufmfam=\teneufm
\scriptfont\eufmfam=\seveneufm
\scriptscriptfont\eufmfam=\fiveeufm

\edef\msafam@{\hexnumber@\msafam}
\edef\msbfam@{\hexnumber@\msbfam}
\def\Bbb#1{\fam\msbfam\relax#1}
\def\frak#1{{\fam\eufmfam\relax#1}}

\def\newsymbol#1#2#3#4#5{\let\next@\relax
 \ifnum#2=\@ne\let\next@\msafam@\else
 \ifnum#2=\tw@\let\next@\msbfam@\fi\fi
 \mathchardef#1="#3\next@#4#5}

\newsymbol\square 1003
\newsymbol\twoheadrightarrow 1310

\catcode`\@=12

\def\qed{\hfill$\qquad\square$\par\medbreak}


\magnification=1200
\hsize=150truemm
\hoffset=5truemm
\vsize=235truemm

\def\x{{\bar x}}

\def\sumj#1{\sum_{i=0}^{#1}}
\def\sumi{\sumj{2g+2}}

\def\auf#1{{\left\lceil #1\over 2 \right\rceil}}
\def\ab#1{{\left\lfloor #1\over 2 \right\rfloor}}

\def\wlx{\overline{X}}
\def\wly{\overline{Y}}

\def\wt #1{\widetilde{#1}}
\def\wl #1{\overline{#1}}
\def\tf{{\tilde f}}
\def\mdot#1#2{#1$, \dots,~$#2}
\def\sier#1{{\cal O}_{#1}}

\def\C{{\Bbb C}}
\def\P{{\Bbb P}}
\def\Z{{\Bbb Z}}

\def\D{{\cal D}{\it iff}}
\def\ttw#1{T^2_X(#1)}
\def\tee#1{T^1_Y(#1)}

\def\coker{\mathop{\rm coker}}
\def\Hom{\mathop{\rm Hom}}
\def\Sym{\mathop{\rm Sym}\nolimits}
\def\Pic{\mathop{\rm Pic}}

\def\Rank{\mathop{\rm Rank}}
\def\elm{\mathop{\rm elm}}
\def\mapdown#1{\Big\downarrow
 \rlap{$\vcenter{\hbox{$\scriptstyle #1$}}$}}
\def\mapright#1{\smash{\mathop{\longrightarrow}\limits^{#1}}}

\def\cdmatrix#1{\def\normalbaselines{\baselineskip20pt
 \lineskip3pt \lineskiplimit3pt }
\matrix{#1}}

\def\cite#1{[#1]}
\def\roep(#1)#2.{\medbreak\noindent(#1)%
\enspace{\sl #2\/}.\enspace}
\def\bewijs{\medbreak\noindent{\sl Proof\/}.\enspace}
\def\endroep{\par\medbreak}

\def\beginthmenv#1#2(#3){\medbreak\noindent
(#3)\enspace{\bf #1#2.\enspace}\begingroup\sl}
\def\endthmenv{\endgroup\par\medbreak}

\def\beginprop{\beginthmenv{Proposition}}
\def\endprop{\endthmenv}
\def\beginlem{\beginthmenv{Lemma}}
\def\endlem{\endthmenv}
\def\beginthm{\beginthmenv{Theorem}}
\def\endthm{\endthmenv}
\def\begincor{\beginthmenv{Corollary}}
\def\endcor{\endthmenv}

\outer\def\beginsection#1\par{\bigskip\bigskip
\message{#1}\leftline{\largebf#1}\nobreak
\medskip\noindent}

\def\rfac{Arbarello--Cornalba 1981}
\def\rfas{Arbarello--Sernesi 1978}
\def\rfmac{Bayer--Stillman}
\def\rfbc{Behnke--Christophersen 1991}
\def\rfcas{Castelnuovo 1890}
\def\rfchm{Ciliberto--Harris--Miranda 1988}
\def\rfco{Conforto 1939}
\def\rfdr{Drewes 1989}
\def\EGA{EGA IV}
\def\rfgr{Green 1984}
\def\rfgm{Greuel 1982}
\def\rfHAG{Hartshorne 1977}
\def\rfkl{Kleiman 1977}
\def\rflev{L\'evy-Bruhl-Mathieu 1953}
\def\rflo{Looijenga 1984}
\def\rflw{Looijenga--Wahl 1986}
\def\rfmum{Mumford 1973}
\def\rfTata{Mumford 1984}
\def\rfpie{Piene 1978}
\def\rfpi{Pinkham 1970}
\def\rfr{Reid 1989}
\def\rfschl{Schles\-singer 1973}
\def\rfschr{Schreyer 1986}
\def\rfse{Severi 1921}
\def\rfstb{Stevens 1991}
\def\rfstp{Stevens 1989}
\def\rftend{Tendian 1992a}
\def\rftenh{Tendian 1992b}
\def\rftenw{Tendian 1990}
\def\rfwad{Wahl 1988}
\def\rfwaj{Wahl 1987}
\def\rfwag{Wahl 1990}
\def\rfwai{Wahl 1989}


\vglue 1cm
\centerline{\LARGEbf
Deformations of cones}
\medskip
\centerline{\LARGEbf
over hyperelliptic curves}
\bigskip
\centerline{\largebf Jan Stevens}

\beginsection

The deformation theory of a two-dimensional singularity, which is
isomorphic to an affine cone over a curve, is intimately linked with the
(extrinsic) geometry of this curve. In recent times various authors have
studied one-parameter deformations, partly under the guise of extensions
of curves to surfaces (cf.~the survey \cite{\rfwai}).  In this paper we
consider the versal deformation of cones, in the simplest case: cones
over hyperelliptic curves of high degree. In particular, we show that
for  degree $4g+4$, the highest degree for which interesting deformations
exist, the number of smoothing components is $2^{2g+1}$ (the case $g=3$
is exceptional).

Let $X$ be the cone over a hyperelliptic curve $C$, embedded with a line
bundle $L$. If $d:=\deg L \geq 2g+3$, then $X$ has always infinitesimal
deformations in negative degrees:   $\dim T^1_X(-1)=2g+2$ \cite{\rfdr}.
On the other hand, one knows that only conical deformations exist  over
reduced base spaces, if $\deg L>4g+4$ \cite{\rftend}. This implies that
all deformations in negative degree must be obstructed. If $S$ is a
surface with $C$ as hyperplane section, then one can degenerate $S$ to the
projective cone over $C$, or from another point of view, deform the
projective cone over $C$ to $S$; Pinkham calls this construction
`sweeping out the cone' \cite{\rfpi}. Surfaces with hyperelliptic
hyperplane sections were already classified by Castelnuovo, and the
supernormal surfaces among  them have degree $4g+4$ \cite{\rfcas}. They
are rational ruled surfaces, and such surfaces come in two deformation
types; therefore there are at least two smoothing components. This
observation was the starting point of the present paper. A computer
computation of the versal deformation in negative degree  with {\sl
Macaulay\/} \cite{\rfmac} gave for  an example with $g=2$ the number of
$32$ smoothing components.

As  the versal deformation can be chosen  $\C^*$-equivariant, it makes
sense to restrict to the part of negative degree. We want to show that
the base space $S^-$ has $2^{2g+1}$ one-dimensional components. First of
all we have to exhibit this number of surfaces with $C$ as hyperplane
section, such that  the normal bundle of $C$ in the surface is $L$. The
main point is that an elementary transformation on the ruled surface $S$
in a Weierstra\ss\ point of $C$ does not change the normal bundle of $C$.
The composition of elementary transformations in all  Weierstra\ss\
points gives an involution on $S$; we get $2^{2g+2}/2$ surfaces. This
construction works for every hyperelliptic curve $C$ and every line
bundle $L$ on $C$. We obtain $2^{2g+1}$ smooth subspaces of dimension
$3g$ of the base of the versal deformation; by a result of \cite{\rftend}
this dimension is exactly the dimension of smoothing components.
Therefore we have  found $2^{2g+1}$ irreducible components.

The next thing to determine is $T^2$, the space where the obstructions
lie. For $T^2(\nu)$, $\nu<-2$, we have a general vanishing result
\cite{\rfwaj}, and it is not difficult to find the dimension in case
$\nu>-2$. We find the dimension of $T^2$ with the {\sc Main Lemma} of
\cite{\rfbc}, which connects the number of generators of $T^2$ with the
codimension of smoothing components in the base space of the versal
deformation of a general hyperplane section.  Therefore we compute the
dimension of $T^1$ for the cone over $d$ points on a rational normal
curve of degree $d-g-1$ in $\P^{d-g-1}$. We use explicit equations for
the curve.  We also need the equations of $X$ to show that the $\sier
X$-module $T^2$ is annihilated by the maximal ideal. Altogether we obtain
that $\dim T^2(-2)=d-2g-3$, $\dim T^2(-1)=(g-2)(d-g-3)$, and the other
$T^2(\nu)$ vanish, if $d>2g+3$.

Actually, the equations for the cone over a hyperelliptic curve, or for
its hyperplane sections, have a nice structure, which Miles Reid calls
the {\sl rolling factors format\/} \cite{\rfr}. We give an interpretation
of $T^2(-2)$ in terms of the rolling factors format, in (2.12). It is not
difficult to compute the part  of the versal deformation in negative
degrees, once one has represented $T^1(-1)$ as perturbations of the
equations. Unfortunately, this representation is given by complicated
formula's, which we only computed in the simplest case, that $L$ is a
multiple of the $g^1_2$. The resulting equations for the base space $S^-$
are as complicated, and it is difficult to see if they define a complete
intersection; for $d=4g+4$ we have $2g+1$ quadratic equations in $2g+2$
variables. By determining explicitly the elements of $T^1$, induced by
our $2^{2g+1}$ surfaces, we show that we have a complete intersection,
and this allows us finally to conclude that there are no other components
(except when $g=3$).

Powerful methods exist to compute $T^1$ for surface singularities,
without using explicit equations. They are based on Schlessinger's
description of $T^1$ \cite{\rfschl}: suppose $\mathop{\rm depth}_0
X\geq2$ and write $U=X\setminus 0$, then
$T^1_X=\ker\bigl\{H^1(U,\Theta_U)\to H^1(U,\Theta_n|_X)\bigr\}$. In the
special case that $X$ is the cone over a projective variety $Y$, all
sheaves are graded, and the graded parts can be computed on $Y$, see
\cite{\rfschl} and Mumford's `footnote' to it \cite{\rfmum}.  We describe
the situation in terms of the sheaf of differential operators of order
$\leq 1$ \cite{\EGA}. In the case of cones over curves it is advantageous
to dualise. Now the bundle of principal parts comes in, and with it
Wahl's Gaussian map (cf.~ \cite{\rfwai}). For cones over curves of high
degree the computation of $T^1(-1)$ is the most difficult; much of the
work on the Gaussian map is connected with this case, and specially with
vanishing results.   The most complete results on interesting
deformations  are obtained by Sonny Tendian \cite{\rftenw}.

With a trick, which basically is contained in \cite{\rfmum}, one sees
that for non hyperelliptic curves, embedded with a non special line
bundle $L$, the dimension of $T^1(-1)$ of the cone over $C$ is equal to
$h^0(C,N_K\otimes L^{-1})$, where $N_K$ is the normal bundle of $C$ in
its canonical embedding. For low genus this gives quite precise
information, because then the normal bundle $N_K$ is easy to describe.
We take the opportunity to remark that the computations in \cite{\rfstp}
imply that for a general curve of genus $g\geq 3$ the Gaussian map
$\Phi_{K,L}$ is surjective for all line bundles $L$ with $\deg L\geq
2g+11$.    In the hyperelliptic case a variant of the construction yields
easily the dimension of $T^1(-1)$.

The paper \cite{\rflw} introduces  a collection of smoothing data for
surface singularities, which in many cases distinguish between smoothing
components. In our case we determine a subset with $1+2^{2g}$ elements,
onto which the set of smoothing components is mapped surjectively, so
these smoothing data do not suffice to distinguish all components. The
computations are similar to the case of a simple elliptic singularity of
degree $8=4g+4$, where the number of smoothing components is really $1+4$.

The organisation of this paper is as follows: in Section 1 we discuss the
description of $T^1$ for cones, and give precise results in the
hyperelliptic case. We also give a formula for the graded parts of $T^2$.
In Section 2 we describe the equations for hyperelliptic cones. We
compute  $T^1$ for the general hyperplane section, and deduce the
dimension of $T^2$ from it. In the last Section we prove the results on
smoothing components.


\beginsection 1. Cones over curves

(1.1) The basic reference for the deformation theory of cones is a paper
by Schles\-singer \cite{\rfschl} and Mumford's `footnote' to it
\cite{\rfmum}.  We start with the description of $T^1_X$ for a singularity
$(X,0)\subset (\C^n,0)$:
$$
0\longrightarrow \Theta_X \longrightarrow \Theta_{n}|_X
\longrightarrow N_X \longrightarrow T^1_X\longrightarrow 0.
\eqno(*)
$$
Here $N_X=\Hom(I/I^2,\sier X)$ is the normal sheaf; the tangent
sheaf of $X$ is also defined as a dual:
$\Theta_X=\Hom(\Omega_X^1,\sier X)$.

If $Z\subset X$ is a closed subset, containing the singular locus ${\rm
Sing}\, X$ of $X$,  and if $\mathop{\rm depth}_Z X\geq2$, then
$T^1_X=\coker\bigl\{H^0(U,\Theta_n|_X )\to H^0(U,N_U)\bigr\}$, or
$T^1_X=\ker\bigl\{H^1(U,\Theta_U)\to H^1(U, \Theta_n|_X)\bigr\}$, where
$U=X\setminus Z$. In the special case that $X$ is the cone over a
projective variety $Y$ we want to interpret these groups on $Y$. From now
on we work in the algebraic category.

Let $Y$ be a smooth projective variety, and $L$ a very ample line bundle
on $Y$. We set $V=H^0(Y,L)$. Let $\phi_L\colon Y\to \P(V^*)$ be the
embedding of $Y$ with $L$.  Let $X\subset V^*$ be the affine cone over
$\phi_L(Y)$. We will identify $Y$ and $\phi_L(Y)$; we have $L={\cal
O}_Y(1)$. Suppose that $\phi_L(Y)$ is projectively normal, i.e. $X$ is
normal. The smooth space $U:=X-0$ is a $\C^*$-bundle over $Y$. We denote
with the same symbol  $\pi$ the projections  $\pi\colon U\to Y$ and
$\pi\colon V^*-0\to \P(V^*)$. If ${\cal F}$ is a sheaf on $X$ with a
natural $\C^*$-action, then $\pi_*{\cal F}$ decomposes into direct sums of
the eigenspaces for the various characters of $\C^*$. Let ${\cal F}_Y$ be
the sheaf of $\C^*$ invariants. Then:
$$
H^0(X,{\cal F})=H^0(U,{\cal F})=
\bigoplus_{\nu=-\infty}^{\nu=\infty} H^0(Y, {\cal F}_Y(\nu)).
$$
The last equality does not hold in the analytic context; in that case one
has a filtration on $H^0(X,{\cal F})$, whose associated graded space is
the direct sum as above, which would suffice for our purposes.

We describe the sheafs in the exact sequence ($*$). Actually, when
explicitly writing down sections, it is more convenient to use
homogeneous coordinates, i.e.~to compute on $X$. We illustrate this with
the tangent bundle of $\P^n$. One has the Euler sequence \cite{\rfHAG,
II.8.20.1}:
$$
0\longrightarrow \sier{\P^n} \longrightarrow
\sier{\P^n}(1)^{n+1} \longrightarrow \Theta_{\P^n}
\longrightarrow 0.
$$
In homogeneous coordinates $(\mdot{z_0}{z_n})$ elements of a basis of
$H^0(\sier{\P^n}(1)^{n+1})$ can be written as $z_i{\partial\over\partial
z_j}$, and the map from $\sier{\P^n}$  is given by $1\mapsto \sum
z_i{\partial\over\partial z_i}$.

We recall the definition of the {\sl sheaf of principal parts\/}
\cite{\rfkl, IV.A, \EGA.16.7}. Consider a scheme $Y$; let ${\cal J}$ be
the ideal sheaf of the diagonal $\Delta$ in $Y\times Y$, and let
$Y_\Delta^{(n)}$ be the n-th infinitesimal neighbourhood of $\Delta$. The
canonical projections $p_1$ and $p_2$ of the product induce maps
$p_1^{(n)}\colon Y_\Delta^{(n)}\to Y$ and $p_2^{(n)}\colon
Y_\Delta^{(n)}\to Y$. We define for a sheaf ${\cal F}$ \cite{\EGA.
16.7.2}:
$$
{\cal P}^n_Y({\cal F})=(p_1^{(n)})_*((p_1^{(n)})^*({\cal F})).
$$
We write ${\cal P}^n_Y$ for ${\cal P}^n_Y(\sier Y)$. One has  ${\cal
P}^n_Y({\cal F})={\cal P}^n_Y\otimes_{\sier Y}{\cal F}$, where the tensor
product is taken with the $\sier Y$-module structure, defined by $p_2$.
Because the diagonal is a section of $Y\times  Y$ for both $p_1$ and
$p_2$, both morphisms define a homomorphism $\sier Y \to {\cal P}^n_Y$,
and therefore an $\sier Y$-module structure.  Except when explicitly
stated, we always consider  the $\sier Y$-module structure on ${\cal
P}^n_Y$, induced by $p_1$, and write it as left multiplication. One
denotes by $d^n$ the morphism $\sier Y \to {\cal P}^n_Y$, induced by
$p_2$ \cite{\EGA.16.3.6}. For every $t\in \Gamma(U,\sier Y)$, $U\subset
Y$ open, $d^nt$ is the {\sl principal part of order $n$\/}. In
particular, $dt=d^1t-t\in \Gamma(U,\Omega^1_Y)$ is the differential of
$t$.  We have the exact sequence:
$$
0\longrightarrow \Omega^1_Y \longrightarrow
{\cal P}^1_Y \longrightarrow \sier Y\longrightarrow 0.
$$
On $\P^n$ this sequence is the dual of the Euler sequence.

The sheaf ${\cal P}^n_Y({\cal F})$  has also two $\sier Y$-module
structures; it is convenient  to write them on the left and the right.
For $a\in \Gamma(U,\sier Y)$, $b\in \Gamma(U,{\cal P}^n_Y)$ and $t\in
\Gamma(U,{\cal F})$ one has \cite{\EGA.16.7.4}:
$$
a(b\otimes t)=(ab)\otimes t, \qquad  (b\otimes t)a=(b\cdot d^na)\otimes
t=b\otimes (at)=(d^na)\cdot(b\otimes t).
$$
There is a map $d^n_{\cal F}\colon {\cal F}\to{\cal P}^n_Y({\cal
F})$ with $d^n_{\cal F}(t)=1\otimes t$.

\roep(1.2)Definition {\rm \cite{\EGA.16.8.1}}. Let ${\cal F}$ and ${\cal
G}$ be two $\sier Y$-modules. A homomorphism $D\colon {\cal F} \to {\cal
G}$ is a {\sl differential operator of order $\leq n$\/} if there exists
a homomorphism  $u\colon{\cal P}^n_Y({\cal F}) \to {\cal G}$ such that $D=
u\circ d^n_{\cal F}$.
\endroep

The differential operators form a group; by applying the construction on
open sets we obtain a sheaf $\D_Y^n({\cal F},{\cal G})$, which is
isomorphic to ${\cal H}{\it om}_{\sier Y} ({\cal P}^n_Y({\cal F}),{\cal
G})$ \cite{\EGA.16.8.4}. We write simply $\D_Y^n$ for $\D_Y^n(\sier
Y,\sier Y)$, and  $\D_Y^n({\cal F})$ for $\D_Y^n({\cal F},\sier Y)$.

We return to our embedding $\phi_L\colon Y\to \P(V^*)$. Then
$\phi_L^*({\cal P}^1_{\P(V^*)})=V\otimes_\C L^{-1}$. Here $V$ is
considered to be the vector space with $dz_i$ as basis. Let $N^*_Y$ be
the conormal bundle of $Y$ in $\P(V^*)$. Because $Y$ is smooth, one has
the familiar exact sequence $0\longrightarrow N^*_Y\longrightarrow
\phi_L^* (\Omega^1_{\P(V^*)}) \longrightarrow \Omega^1_Y \longrightarrow
0$, with which we obtain the following result \cite{\rfkl, (IV.19)}:

\beginprop(1.3)
In the situation as above the following sequence is exact:
$$
0\longrightarrow N^*_Y\otimes L \longrightarrow
V\otimes_\C \sier Y
\longrightarrow {\cal P}^1_Y(L)\longrightarrow 0,
$$
or dually:
$$
0\longrightarrow \D^1_Y(L) \longrightarrow
V^*\otimes_\C \sier Y
\longrightarrow N_Y\otimes L^{-1}\longrightarrow 0.
$$
\endprop

In particular we can view the $\C^*$-invariants of the exact sequence
($*$) as obtained by taking global sections of the sequence
$0\longrightarrow \D^1_Y \longrightarrow  V^*\otimes_\C L \longrightarrow
N_Y \longrightarrow 0$ on $Y$. We get the following formulation of a
result  of \cite{\rfschl}:

\beginthm(1.4)
Let $L$ be a very ample line bundle on  a smooth projective
variety  $Y$, which embeds $Y$ as projectively normal subvariety
of\/ $\P(V^*)$, where $V=H^0(Y,L)$.  Let $X\subset V^*$ be the
affine cone over $Y$.  Then the graded parts $T^1_X(\nu)$ of
$T_X^1$ are given by:
$$
T^1_X(\nu)= \coker \{V^*\otimes H^0(Y,L^{\nu+1})
\longrightarrow H^0(Y,N_Y\otimes L^\nu) \},
$$
or alternatively,
$$
T^1_X(\nu)= \ker \{H^1(Y,\D_Y^1\otimes L^\nu)
\longrightarrow V^*\otimes H^1(Y,L^{\nu+1}) \}.
$$
\endthm

{}From now one we concentrate on  the case that $Y$ is one dimensional;
Serre duality then transforms the second formula into one involving $H^0$.

\begincor(1.5)
Let  $X$ be the cone over a curve $C$, embedded by $L=\sier C(1)$.
Write $K$ for $\Omega^1_C$. Then:
$$
(T^1_X)^*(\nu)= \coker \{V\otimes H^0(C,K(-\nu-1))
\longrightarrow H^0(C,{\cal P}^1_C\otimes K(-\nu))
\}.
$$
\endcor

This Corollary  makes it possible to determine the graded parts of $T^1$
in many cases. Vanishing results exist for line bundles of high degree.
We recall some results for the various degrees $\nu$. We assume that
$g(C)\geq 2$.

\roep (1.6)Case\/ {\rm I:} $\nu\geq1$. Suppose first that $L^\nu$ is non
special, i.e.~$H^1(C,L^\nu)=0$.
Then   $T^1_X(\nu)= H^1(C,\D_C^1\otimes L^\nu)
=H^1(C,\Theta_C\otimes L^\nu)$. In particular, if $\deg L > 4g-4$,
then  $T^1_X(\nu)=0$ \cite{\rfmum}. If $L^\nu$ is special, but
$H^1(C,L^{\nu+1})=0$, then still $T^1_X(\nu)= H^1(C,\D_C^1\otimes
L^\nu)$.  One has the exact sequence:
$$
0\longrightarrow
H^1(C,L^\nu) \longrightarrow
H^1(C,\D_C^1\otimes L^\nu) \longrightarrow
H^1(C,\Theta_C\otimes L^\nu)\longrightarrow0.
$$
If $L^\nu=K$, this follows  from the surjectivity of $\Gamma({\cal
P}^1_C) \to \Gamma(\sier C)$.

\roep (1.7)Case\/ {\rm II:} $\nu=0$.
One has  $T^1_X(0)= \ker \{H^1(C,\D_C^1) \longrightarrow V^*\otimes
H^1(C,L) \}$. The map here can also be thought of as a cup product
$H^1(C,\D_C^1)\otimes H^0(C,L)\longrightarrow H^1(C,L) $ \cite{\rfac};
for a differential operator $\theta$ and a section  $s\in \Gamma(L)$ the
cup product $\theta\cdot s=0\in H^1(L)$ if and only if the section $s$
lifts to the first order deformation of $L\to C$, defined by $\theta$. In
\cite{loc.cit.} the vector space $T^1_X(0)$ is identified as tangent
space to a space of $g^r_d$'s on a variable curve. More precisely, let
$\pi\colon X\to S$ be a miniversal family of smooth curves, and consider
the relative Picard variety $\Pic_{X/S}^d$ and the bundle ${\cal G}^r_d$
over $S$ of  $g^r_d$'s on the fibres; there is a map $c\colon {\cal
W}^r_d\to \Pic_{X/S}^d$ with image  ${\cal W}^r_d$. Because by assumption
the linear series $L$ is complete, the map $c$ is injective. Then $H^1(C,
\D_C^1)$ is the tangent space to $\Pic_{X/S}^d$ in the point $L\to C$,
and  $T^1_X(0)$ is the tangent space to ${\cal G}^r_d$.  The problem now
is to prove that ${\cal G}^r_d$ is smooth of expected dimension
$3g-3+\rho$.  However, for $\rho<0$ not much is known.

For non special $L$ one has $\dim T^1_X(0)=4g+3$. In terms of
deformations this means that the cone can be deformed by changing the
moduli of the curve or by changing the line bundle $L$ in $\Pic^d$.  For
$L=K$ the composed map $H^1(\sier C)\to H^1(\D_C^1)\to V^*\otimes
H^1(K)\cong H^1(\sier C)$ is an isomorphism, so $\dim T^1_X(0) = 3g-3$
and the versal family is a family of curves in their canonical embedding.

\roep (1.8)Case\/ {\rm III:} $\nu\leq-2$.
In this case $T^1_X(\nu)=H^0(C,N_C \otimes L^{\nu})$ because $L^{\nu+1}$
is a line bundle of negative degree. We have two vanishing results.

\beginlem {\rm \ \cite{\rfmum}}(1.9)
If\/ $T^1_X(-1)=0$, then $T^1_X(\nu)=0$ for all $\nu\leq-2$.
\endlem

\bewijs
If $\Gamma(N_C\otimes L^\nu)\neq0$, then $N_C\otimes L^{-2}$ has a
non-zero section $s$, and for all $t\in V=\Gamma(L)$ the tensor product
$t\otimes s$ is a non-zero section of $N_C\otimes L^{-1}$; therefore
$h^0(C,N_C\otimes L^{-1})\geq \dim V$. Because $T^1_X(-1)=0$, the map
$V^*\to\Gamma(N_C\otimes L^{-1})$ is surjective, so all sections are of
the form $t\otimes s$, therefore they are proportional, and do not
generate $N_C\otimes L^{-1}$. But $N_C\otimes L^{-1}$ is generated by its
sections, because $V^*\otimes \sier C$ is, and $V^*\otimes \sier C\to
N_C\otimes L^{-1}$ is surjective.
\qed

\roep (1.10)Remark. The argument in the proof shows that
$\Gamma(N_C\otimes L^\nu)=0$, $\nu\leq-2$, if $\dim T^1_X(-1)< {\rm
rank}\,N_C-1$.

\beginprop {\rm \ \cite{\rfwaj, 2.5}}(1.11)
Let $Y\subset\P=\P(V^*)$, with $V\subset \Gamma(L)$, be a projective
variety, defined by a system of quadratic equations $\boldmath f$.
Suppose every non-zero quadratic equation $f_i$ is involved in a linear
relation $r$; this is true if the  relations are generated by linear ones.
Then $H^0(Y,N_Y \otimes L^{\nu})=0$ for $\nu\leq-2$.
\endprop

\bewijs
Let ${\cal I}$ be the ideal sheaf of $Y$. Consider the complex:
$$
\sier\P(-3)^{\oplus l}\mapright{\boldmath r}\sier\P(-2)^{\oplus
k} \mapright{\boldmath f} {\cal I}\longrightarrow 0,
$$
which  is not necessarily exact at $\sier\P(-2)^{\oplus
k}$. Dualise, twist and restrict to $Y$ to get:
$$
0\longrightarrow N_Y\otimes L^{-2} \longrightarrow
\sier Y^{\oplus k} \mapright{{}^t\boldmath r}
\sier Y(1)^{\oplus l}.
$$
Let $K$ be the kernel of the map ${}^t{\boldmath r}\colon\C^k \to
\Gamma(L)^{\oplus l}$. Here $\C^k$ can be identified with  the dual of
the vector space $Q$ of quadratic equations. The relations $\boldmath r$
involve only equations in $K^\perp$. Therefore $K=0$.
\qed

By a theorem of Green  the conditions are satisfied for an embedding of a
curve with a complete linear system of degree $d\geq 2g+3$ \cite{\rfgr,
Thm.~4.a.1}.

\roep (1.12)Example.
Let $C$ be hyperelliptic, with involution $\pi\colon C\to\P^1$, and let
$L$ be very ample of degree $d$. Then $\phi_L(C)$ lies on a scroll $\wl
S$, the image in $\P(\Gamma(C,L)^*)$ of $S=\P(\pi_* L)$, where $\pi_*
L\cong \sier {}(a)\oplus\sier {}(b)$ with $a+b=d-(g+1)$, and $a,b\leq
d/2$. Suppose $b\leq a$, write $e=a-b$, so $0\leq e\leq g+1$, then
$S\cong\P(\sier{}\oplus\sier{}(-e))$. The Picard group of $S$ is
generated by the section $E_0$ with $E_0^2=-e$, and the class $f$ of a
fibre. We have $C\sim 2E_0+(g+1+e)f$ (the coefficient of $f$ can be
computed from the adjunction formula). Therefore $C^2=4g+4$.

Now suppose that $d=\deg L=2g+2$. Then $\wl S=S$, except when $b=0$,
which occurs if $L=(g+1)g^1_2$; in that case $X$ is the cone over a
rational normal curve of degree $g+1$, and $C$ does not pass through the
vertex (because $E_0\cdot C=0$). We have the normal bundle exact sequence
$0\to N_{C/S} \to N_C \to N_{\wl S}|_C \to 0$. Because the scroll $\wl S$
is defined by quadratic equations with linear relations, the argument of
the proposition gives that $H^0(C, N_{\wl S}\otimes L^{-2})=0$. Therefore
$\Gamma(N_C(-2))=\Gamma(N_{C/S}(-2))$. Because $L\sim E_0+af$, we have
$C\cdot (C-2L)=C\cdot (g+1+e-2a)f=0$, so $N_{C/S}(-2) \cong \sier C$ and
$h^0(C,N_{C/S}(-2))=1$. For $L=(g+1)g^1_2$ this result was obtained by
Drewes \cite{\rfdr}.

We specialise to the case $g=3$. We proved in \cite{\rfstb}, that the
cone  $X$  over $\phi_L(C)$ has 3 smoothing components, if
$L=2K(=4g^1_2)$, and 2 components otherwise. The curve $C$ is a complete
intersection of the scroll and a quadric. If $L=2K$, the scroll is the
projective cone over a rational normal curve of degree 4, which itself
has two smoothing components; the deformation to the Veronese surface
occurs, if we deform $C$ to a non hyperelliptic curve, with $L=2K$.

We compute $T^1_X(-2)$ for $g(C)=3$, $\deg L=8$, $C$ not hyperelliptic.
If $L=2K$, we find as above that
$\Gamma(N_C(-2))=\Gamma(N_{C/S}(-2))=\Gamma(\sier C)$, where $S$ is the
Veronese surface. If $L\neq 2K$, let $D$ be a general divisor in the
linear system; it  is cut out on the canonical curve $C_4$  in $\P^2$ by
a cubic $C_3$, and the linear system is the system of cubics through the
residual intersection $C_4\cap C_3 -D$. Therefore $C\subset \P^5$ lies on
(non unique) Del Pezzo surface of degree 5. One checks that all equations
(which are quadratic) occur in linear relations, although the relations
are not generated by linear ones. So $T^1_X(-2)=0$. A smoothing is
obtained by sweeping out the cone over the blow-up of $\P^2$ with the
linear system of quartics with base points in $D$.  By \cite{\rftend,
Prop. 6.1} the dimension of the corresponding smoothing component  is 15,
which is also the dimension of $T^1$ (see the description of $T^1(-1)$
below), so the base space is smooth. For $L=2K$ this component has
codimension 1, and there is a second component; if $C=S\cap\{ Q=0\}$, then
$S\cap\{ Q=t^2\}$ is $\P^2$, branched in $C_4$.

\roep (1.13)Case\/ {\rm IV:} $ \nu=-1$.
This is the most difficult case.
A more specific knowledge of the maps in the theorem
is necessary. We review the relation with
Wahl's Gaussian map \cite{\rfwag}, see also \cite{\rftenw, \rfdr, \rfwai}.

We start with the following diagram:
$$
\cdmatrix{
&&&&V\otimes\sier C&&&&\cr
&&&&\Big\downarrow&&&&\cr
0& \to & K\otimes L & \longrightarrow & {\cal P}^1_C(L) &
\longrightarrow & L & \to & 0 \cr
}
$$
The kernel of the composed map $V\otimes\sier C\to L$ is a vector bundle
${\cal M}_L$ over $C$, and we get a map  ${\cal M}_L\to K\otimes L$. Let
$M$ be a second line bundle, and tensor everything with $M$. Define:
$$
{\cal R}(L,M)=\ker \{ \mu_{L,M}\colon \Gamma(L)\otimes \Gamma(M)\to
\Gamma(L\otimes M) \}.
$$
Then  ${\cal R}(L,M)=\Gamma({\cal M}_L\otimes M)$ and we have the {\sl
Gaussian map\/} $\Phi_{L,M}\colon {\cal R}(L,M) \to \Gamma(K\otimes
L\otimes M)$. This map is given explicitly by $\Phi_{L,M} (\alpha)=\sum
d^1_Ll_i\otimes m_i$, where $\alpha=\sum l_i\otimes m_i\in{\cal R}(L,M)$
with  $l_i\in \Gamma(L)$ and $m_i\in \Gamma(M)$. A more symmetric
definition can be given in local coordinates; represent sections on an
open set $U$ by functions, again denoted by $l_i$ and $m_i$. Then from
$\sum l_im_i=0$ we get $\sum (l_i\,dm_i+m_i\,dl_i)=0$, so $\Phi_{L,M}
(\alpha)$ can be represented by the $1$-form
$1/2\sum(l_i\,dm_i-m_i\,dl_i)$.

For curves we have that $H^1(K\otimes L\otimes M)=0$, so if we denote the
map $\Gamma(L)\otimes \Gamma(M)\to \Gamma({\cal P}^1_C (L) \otimes M)$ by
$d^1_L\otimes 1_M$, we have the exact sequence $0\to \coker \Phi_{L,M}\to
\coker d^1_L\otimes 1_M \to \coker \mu_{L,M} \to 0$.

We can also start the construction above with $M$ in stead of $L$.  Up to
a permutation of  factors one has $\mu_{L,M}=\mu_{M,L}$ and
$\Phi_{L,M}=-\Phi_{M,L}$. Therefore:
$$
\coker d^1_L\otimes1_M\cong \coker d^1_M\otimes 1_L.
$$

We apply this to the computation of $T^1_X(-1)$. From
Corollary (1.5) we have $(T^1_X(\nu))^*=\coker d^1_L\otimes
1_{K\otimes L^{-(\nu+1)}}$ and in particular $(T^1_X(-1))^*=\coker
d^1_L\otimes 1_K$. Therefore $(T^1_X(-1))^*\cong\coker
d^1_K\otimes 1_L$.

If $\phi_L\colon C\to \P(H^0(C,L)^*)$ is birational onto its image, then
the multiplication maps $\mu_{L,K\otimes L^{-(\nu+1)}}$ are surjective
\cite{\rfas, Proof of Thm 1.6}. This applies to our situation, so
$(T^1_X(\nu))^*=\coker \Phi_{L,K\otimes L^{-(\nu+1)}}$ and in particular
$(T^1_X(-1))^*=\coker\Phi_{L,K}=\coker\Phi_{K,L}$.

Now suppose that $C$ is {\sl not hyperelliptic\/}. Then $K$ is
very ample and we have the exact sequence:
$$
0\longrightarrow N^*_K\otimes K \longrightarrow
\Gamma(K)\otimes_\C \sier C
\longrightarrow {\cal P}^1_C(K)\longrightarrow 0,
$$
where $N^*_K$ is the conormal bundle of the canonical embedding. This
yields that $(T^1_X(-1))^*=\ker \{ H^1(N^*_K\otimes K \otimes L)\to
\Gamma(K)\otimes H^1(L)\}$. This map is surjective, because $H^1( {\cal
P}^1_C(K)\otimes L)=0$. Dually,  $T^1_X(-1)=\coker \{\Gamma(K)^*\otimes
H^0(K\otimes L^{-1})\to H^0(N_K\otimes L^{-1})\}$. In particular:

\beginlem(1.14)
For the cone $X$ over a non hyperelliptic curve, embedded by a non special
complete linear system,
$T^1_X(-1)=H^0(N_K\otimes L^{-1})$.
\endlem
As Mumford remarks, there is an integer $d_0$, depending only on
$C$, such that $H^0(N_K\otimes L^{-1})=0$ for $\deg L\geq d_0$.
The bundle $\cal E$ of \cite{\rfmum} is $N_K\otimes
K^{-1}$; in fact, Mumford's construction can be understood as
interchanching $K$ and $L$.

In general it is difficult to give sharp explicit bounds, but for low
genus the result is very effective. Before we give some examples, we
prove a  result \cite{\rfwai, 2.2}, which is based on a more general
lemma of Lazarsfeld, giving sufficient conditions for surjectivity of
Gaussian maps.

\beginlem(1.15)
Suppose $C$ is not hyperelliptic, trigonal or a plane quintic. If
$H^0(K^2\otimes L^{-1})=0$ {\rm (}in particular, if $\deg L >4g-4${\rm )},
then $\Phi_{K,L}$ is surjective.
\endlem

\bewijs
By Petri's Theorem the ideal $I$ of the canonical curve $C$ is generated
by quadrics, so there is a surjection $\sier {}(-2)^{\otimes a} \to I$.
Dualising, restricting to $C$ and twisting gives an injection
$H^0(N_K\otimes L^{-1})\to (H^0(K^2\otimes L^{-1}))^{\otimes a}$.   \qed

We now look at low genera.

\roep (1.16)$g=3$. The normal bundle is the line bundle $K^{\otimes4}$.
Therefore $d_0=17$. For $\deg L=16$ only $L=K^{\otimes4}$ has non zero
$T^1(-1)$, and the cone $X(C, K^{\otimes4})$ is smoothable, as hyperplane
section of $X(\P^2, \sier {}(4))$. For all $L$ with $\deg L=13$ the
Gaussian map $\Phi_{K,L}$ is not surjective, for the general $L$ of
degree 14 it is.

\roep (1.17)$g=4$. The normal bundle is $K^{\otimes2}\oplus
K^{\otimes3}$, so $d_0=19$.

More generally, for $C$ trigonal, the canonical curve sits on a scroll
$S$  as divisor of type $3H-(g-4)R$ \cite{\rfschr}, and we have the
normal bundle sequence  $0\to N_{C/S}\to N_C \to N_S|_C\to0$. Here
$N_{C/S}=K^3\otimes (g^1_3)^{4-g}$, a bundle of degree $3g+6$. Sonny
Tendian proves that $H^1(N_S|_C\otimes L^{-1})=0$, if $H^0(K^2\otimes
L^{-1})=0$, or for $L=K^2$ \cite{\rftenw}. In particular, $d_0\leq
\max(3g+7,4g-4)$. We conjecture that $H^1(N_S|_C\otimes L^{-1})=0$ if
$\deg L \geq d_1$ for some $d_1<3g+3$. For $\deg L =3g+6$ we would then
have that $T^1(-1)\neq0$ only if $L$ is the characteristic linear system
of the family of curves of type $3H-(g-4)R$ on $S$, which also gives a one
parameter smoothing.

\roep (1.18)$g=5$. The general curve is a complete intersection of three
quadrics, with normal bundle $3K^2$, and $d_0=17$ (for a trigonal ${}^5C$
we have $d_0=22$).

\roep (1.19)$g=6$. Here the famous possibility of a plane quintic occurs.
It gives $d_0=26$. The general curve lies on a (possibly singular) Del
Pezzo surface $S$ of degree 5 with normal bundle $N_{C/S}=K^2$.

Tendian has studied the case $L=K^2$  \cite{\rftenw}. He proves that
$\Phi_{K,K^2}$ is surjective, if the Clifford index of $C$ is at least 3.
This leaves the hyperelliptic, trigonal and  tetragonal curves, as well
as smooth plane quintics and sextics. We now concentrate on the tetragonal
case. Then $C$ lies on a three dimensional scroll $X$ of type
$S(e_1,e_2,e_3)$ of degree $e_1+e-2+e_3=g-3$ \cite{\rfschr}, but more
important are the surfaces on which $C$ lies, cf. \cite{\rflev}. We
obtain $C$ as complete intersection of divisors $Y\sim 2H-b_1R$ and
$Z\sim 2H-b_2R$ on $X$ with $b_1+b_2=g-5$ and $b_2\leq b_1$. Either $Y$
is rational, or the given $g^1_4$ is composed with an elliptic or
hyperelliptic involution $C\mapright{2:1}E\mapright{2:1}\P^1$, and $Y$ is
a ruled surface over $E$, with a rational curve of double points, which
is the canonical image of $E$; then $\deg Y=g-1+b_2$,
$p_a(E)=1/2(b_2+2)$, and $C$ does not intersect the double curve. The
normal bundle of $C$ in $Y$ is $N_{C/Y}=K^2\otimes (g^1_4)^{-b_2}$. In
the case $b_2=0$, and $g>5$, $c$ lies on a unique elliptic cone $Y$, and
$H^1(C,N_Y\otimes K^{-2})=0$ \cite{\rftenw, 2.2.12}, so $\dim
T^1_{X(C,K^2)}(-1)=1$ and $Y$ leads to a one parameter deformation to a
simple elliptic singularity of degree $g-1$.

\medskip

The examples above, where $\Phi_{K,L}$ is not surjective, have in common
that $N_K$ is unstable; the rank of $N_K$ is $g-2$, the degree is
$2(g^2-1)=(g-2)(2g+4)+6$, and we found surfaces $Y$ with $\deg
N_{C/Y}>2g+4+\lfloor6/g-2\rfloor$. This motivates:

\roep (1.20)Question. Suppose that the canonical curve ${}^gC$ lies on a
unique surface $Y$ with $d_0= \deg N_{C/Y}>2g+4+\lfloor6/g-2\rfloor$. Is
then $\Phi_{K,L}$ surjective for all $L$ with $\deg L>d_0$, and for $L$
with $\deg L=d_0$ and $L\neq N_{C/Y}$?
\endroep

For a general curve such surfaces do not exist; already for
$g=7$ the general curve has a plane representation as
$C^7(8A^2)$, so $\deg N_{C/Y}=17$.

\roep (1.21)Question {\rm \cite{\rfwai, 2.5}}. Is $\Phi_{K,L}$
surjective on a general curve of genus $g\geq 12$ for $\deg L\geq
2g-2$? What is the best bound?
\endroep

We remark that for general $C$ with $g=10$ or $g\geq 12$ the map
$\Phi_{K,K}$ is surjective \cite{\rfchm}. The problem is of course to
find a general curve; the easiest to handle are singular curves. In
\cite{\rfstp} we constructed for all $g\geq5$ a $g$-cuspidal canonically
embedded curve $\Gamma$ with a Weil divisor $D$ of degree $g+5$ such that
$h^0(N_\Gamma(-D))=6$.

\beginprop(1.22)
For a general curve of genus $g\geq 3$ the map $\Phi_{K,L}$
is surjective  for all $L$ with $\deg L\geq 2g+11$.
\endprop

\bewijs
{}From \cite{\rfstp, Prop.~5}  it follows that $h^0(N_C\otimes L^{-1})\leq
6$  for a general curve $C$ with $g\leq5$ and a general $L$ with $\deg L
=g+5$. Therefore $h^0(N_C\otimes L^{-1})=0$, if $\deg L\geq g+11$;
presumably a general $L$ of degree $g+6$ will do for large $g$. Fix a
line bundle $L_0$ of degree $g+11$ with $h^0(N_C\otimes L_0^{-1})=0$.
Every line bundle $L$ with $\deg L\geq 2g+11$ can be written as
$L=L_0(D)$ with $D$ an {\sl effective\/} divisor of degree $\deg L
-g-11\geq q$. This implies that   $H^0(N_C\otimes L^{-1})$ is a subspace
of $H^0(N_C\otimes L_0^{-1})$. For $g=3,4$ see (1.16) and (1.17). \qed

\roep (1.23)The hyperelliptic case.
The formula:
$$
(T^1_X(-1))^*=\coker
d^1_L\otimes 1_K\cong \coker d^1_K\otimes 1_L
$$
still holds.  However the map $d^1_K\colon \Gamma(K)\otimes \sier C\to
{\cal P}^1(K)$ is only generically surjective, with cokernel  equal to
$\coker \{ \phi_K^*\Omega^1_\P\otimes K\to K^2\}\cong \Omega_\phi\otimes
K$, where $\Omega_\phi$ are the relative differentials, so the cokernel
is given by the Jacobian ideal \cite{\rfpie}. Therefore, if $B$ denotes
the set of Weierstra\ss\ points, we have an exact sequence:
$$
\Gamma(K)\otimes \sier C\mapright{d^1_K} {\cal P}^1(K)
\longrightarrow \bigoplus_{p\in B} \C_p \to 0 .
$$
Let ${\cal P}={\rm Im} d^1_K$; this is a rank 2 vector bundle.

\beginlem(1.24)
Let the rational normal curve $R=\phi_K(C)$ be the canonical
image of $C$, and write $\phi\colon C\to R$; let $H=\sier
R(1)\cong \sier {\P^1}(g-1) $. Let $N_R$ be the normal bundle of
$R$ in $\P^{g-1}$. The following sequence is exact:
$$
0\to \phi^*(N^*_R\otimes H^{-1}) \longrightarrow
\Gamma(K)\otimes \sier C\mapright{d^1_K} {\cal P}
\to 0.
$$
\endlem

\bewijs
Consider the commutative diagram:
$$
\cdmatrix{
0&\to &\phi^*(N^*_R\otimes H^{-1})& \longrightarrow
&\Gamma(K)\otimes \phi^*(\sier
R)&\mapright{\phi^*\circ d^1_H} &{\cal P}^1_R(K) & \to &0
\cr &&\Big\downarrow&&\mapdown{\cong}&&\Big\downarrow&&\cr
0&\to& {\cal E}& \longrightarrow &
\Gamma(K)\otimes \sier C&\mapright{d^1_K}& {\cal P}^1_C(K)
&& \cr}
$$
Because $\phi^*\circ d^1_H=d^1_K\circ \phi^*$ \cite{\EGA.16.4.3.4}
the image of the right-hand vertical map is $\cal P$, and therefore
${\cal E}=\phi^*(N^*_R\otimes H^{-1})$.   \qed

The normal bundle of $R$ splits a direct sum of $g-2$ bundles of degree
$g+1$, and therefore ${\cal E}=\oplus_{g-2} (g^1_2)^{-2}$. Using the fact
that $H^1({\cal P}^1_C(K)\otimes L)=0$ for $\deg L>0$, we obtain the
following result; for $\deg L >2g+1$ cf. \cite{\rfwad, 7.11.1, \rfdr}.

\beginprop(1.25)
Let $C$ be a hyperelliptic curve of genus $g$, and $L$ a line
bundle with $\deg L >0$. Then $\dim \coker d^1_K\otimes 1_L=
2g+2+(g-2)h^1(L\otimes(g^1_2)^{-2})-gh^1(L)$. In particular, if
$\deg L>2g+2$, then $\dim T^1(-1)=2g+2$.
\endprop

We remark that $H^1(L)=0$, if $L$ defines a birational map of $C$. Then
$h^1(L\otimes(g^1_2)^{-2})\leq1$, as one sees from the exact sequence
$0\to L\otimes(g^1_2)^{-2} \to L \to \C^4 \to 0$; if
$h^0(L\otimes(g^1_2)^{-2})=0$, this follows because $h^0(L)\geq 3$;
therefore we may assume that
$H^0(L\otimes(g^1_2)^{-2})=\{\mdot{\psi}{\psi z^k}\}$ for some function
$\psi$ on $C$ and a local coordinate $z$ on $\P^1$. Then $H^0(L)$
contains in addition the section $\psi z^{i+1}$, $\psi z^{i+2}$ and at
least one section which is not of this form.

\bigskip

\roep (1.26)$T^2$.
We recall the definition of $T^2_X$ for a singularity $(X,0)\subset
(\C^n,0)$ \cite{\rfschl}: let $0\to  R\to(\sier n)^{\oplus k} \mapright j
I\to 0$ be a resolution of the ideal $I$ of $X$. Let $R_0$ be the
submodule of $R$, generated by the trivial relations, i.e. those of the
form $xj(y)-yj(x)$. Then $R/R_0$ is a $\sier X$-module, and $T^2_X=\coker
\{\Hom((\sier X)^{\oplus k},\sier X)\to \Hom(R/R_0,\sier X)\}$.  We can
also consider the exact sequence $0\to  R_X\to(\sier X)^{\oplus k}  \to
I/I^2\to 0$ on $X$. There  is a surjection $R/R_0\to R_X$, whose kernel
is a torsion module for reduced $X$. Therefore we get the alternative
description:
$$
0\to N_X\longrightarrow \Hom((\sier X)^{\oplus k},\sier X)
\longrightarrow \Hom(R_X,\sier X) \longrightarrow
T^2_X\to 0.
$$
Suppose the sheaf ${\cal T}^2$ has support contained in $Z$ with
$\mathop{\rm dp}_Z X\geq2$, and let  $U=X\setminus Z$.  Then
$T^2_X=\coker\bigl\{H^0(U,\sier n^{\oplus k} )\to H^0(U,R_X^*)\bigr\}
=\ker\bigl\{H^1(U,N_X)\to H^1(U,(\sier n)^{\oplus k}))\bigr\}$.

We specialise to the case of cones as before, so $X$ is the cone over the
smooth projective variety $Y$, embedded with the line bundle $L$. Let the
ideal of $Y$ in $\P(V^*)$, $V=H^0(Y,L)$, be generated by $k$ equations of
degree $\mdot{d_1}{d_k}$. Then the graded parts of $T_X^2$ are given by
the exact sequence:
$$
0\to T^2_X(\nu)\longrightarrow H^1(Y,N_Y(\nu))
\longrightarrow H^1(\sier Y(d_1+\nu))\oplus\cdots\oplus
H^1(\sier Y(d_k+\nu)).
$$
If $\dim Y=1$, then the group $H^1(Y,N_Y(\nu))$ occurs in the
following exact sequence, where we write $C$ for $Y$:
$$
0\to T^1_X(\nu)\longrightarrow H^1(\D_C^1(\nu))
\longrightarrow V\otimes H^1(C,L^{\nu+1})
\longrightarrow  H^1(C,N_C(\nu)) \to 0.
$$
This gives a formula for the dimension of $T^2_X(\nu)$:
$$
\dim T^2_X(\nu)= \dim\ker\{V\otimes H^1(L^{\nu+1})
\to H^1(\oplus_i L^{d_i+\nu})\} -h^1(\D_C^1(\nu))
+\dim T^1_X(\nu).
$$
If $H^1(C,L)=0$, then $H^1(C,N_C(\nu))=0$ for $\nu\geq0$, and therefore
also $T_X^2(\nu)=0$. If $C$ is defined by quadratic equations, and $L$ is
not special, then $T^2_X(-1)=H^1(C,N_C(-1))$ and therefore $\dim
T^2_X(-1)=(g-2)h^0(L)-6(g-1)+\dim T^1_X(-1)$, because $h^1(\D_C^1(-1)) =
2\deg L+4g-4$.

\roep (1.27)Example. The above computation gives new examples of
singularities for which the obstruction map is not surjective, cf.
\cite{\rftenw, 2.4.1}. In particular, if $g(C)=3$ and $L$ is general of
degree $d\geq14$, then $\dim T^1_X(-1)=0$, and $X$ has only conical
deformations, so the base space is smooth, whereas $\dim T^2_X(-1)=d-14$.
For $g=4$ and general $L$ with $9\leq d=\deg L<15$  we have $\dim
T^1_X(-1)=15-d$, $\dim T^1_X=28-d$, which is also the dimension of a
smoothing component  \cite{\rftend, 6.2}. Therefore the base space is
smooth, but $\dim\ttw{-1}=d-9$. Finally, for a general $X({}^5C,L_d)$
with $d>12$ there are only conical deformations, and
$\dim\ttw{-1}=3(d-12)$.
\endroep

For $\ttw\nu$ with $\nu<-2$ we have the following general
vanishing result:
\beginlem {\rm \ \cite{Wahl, Cor. 2.10}}(1.28)
If $Y\subset \P=\P^n$ is a smooth projectively normal subvariety,
defined by quadratic equations, with a resolution
$$
\sier \P(-4)^{\oplus m}\longrightarrow
\sier \P(-3)^{\oplus l}\longrightarrow
\sier \P(-2)^{\oplus k}\longrightarrow
\sier \P\longrightarrow\sier X\longrightarrow 0,
$$
then the cone $X$ over $Y$ satisfies $\ttw\nu=0$ for $\nu<-2$.
\endlem
By \cite{\rfgr, Thm.~4.a.1} the conditions are
satisfied for a curve, embedded with a complete linear system of
degree $d\geq 2g+4$.


\beginsection 2. The versal deformation of hyperelliptic
cones

(2.1) Our aim in this section is to compute equations of the versal
deformation of the cone  $X=X(C,L)$ over a hyperelliptic curve ${}^gC$,
embedded with a complete linear system $L$ of degree at least $2g+3$.
Then, as we have seen, $T^1_X$ is concentrated in degrees $1$, $0$ and
$-1$, and $T^1_X(1)=0$ if $\deg L>4g-4$. We restrict ourselves to the
part of the versal deformation in negative degree, because  otherwise
also non hyperelliptic curves come in.

For the computation of $T^1$  efficient methods exist, which avoid the
explicit  use of equations and relations. For the versal deformation
there seems to be no alternative. Actually, as Frank Schreyer repeatedly
pointed out to me, the equations for $\phi_L(C)$ are rather simple in the
hyperelliptic case: $C$ is a divisor on a two-dimensional scroll, so
besides the determinantal of the scroll we have `essentially' one
equation.

\roep (2.2)Equations for $X$.
Let $L$ be any line bundle of degree $d\geq 2g+3$. Denote the  involution
by $\pi\colon C\to\P^1$. Then $\phi_L(C)$ lies on the scroll
$S=\P_{\P^1}(\pi_* L)$, where $\pi_* L\cong \sier {}(a)\oplus\sier {}(b)$
with $a+b=d-(g+1)$, and $a,b\leq d/2$; in particular, if $L=kg^1_2$, then
$a=k$ and $b=k-(g+1)$. Suppose $b\leq a$, write $e=a-b$, so $0\leq e\leq
g+1$, and $S\cong\P(\sier{}\oplus\sier{}(-e))$. In $\Pic S={\Bbb
Z}E_0\oplus{\Bbb Z}f$, where  $E_0$ is the section with $E_0^2=-e$ and
$f$ is the class of a fibre, we have  $C\sim 2E_0+(g+1+e)f$.

Let $L=kg^1_2$ with $k>g+1$.  Let $C$ be the curve $y^2-\sumi a_ix^i=0$.
We can make this equation quasi-homogeneous by introducing homogeneous
coordinates $(x,\x)$ on $\P^1$, so $C=\{ y^2-\sumi a_i\x^{2g+2-i}x^i=0\}$.
A basis for  $H^0(C,L)$ is given by the functions:
$$
z_i=\x^{k-i}x^i, \quad \hbox{$i=\mdot0k$;}\qquad
w_i=y\x^{k-g-1-i}x^i, \quad \hbox{$i=\mdot0{k-g-1}$.}
$$
The equations for $X$ are those for the scroll, and $2k-2g-1$ further
equations $\phi_m$, obtained by multiplying the equation of $C$ by
suitable powers of $\x$ and $x$, such that the result can be expressed in
the $z_i$ and $w_i$.  These equations are obtained from one another by
{\sl rolling factors\/}, as Miles Reid puts it \cite{\rfr}.  For each
monomial in an equation we get the corresponding monomial in the next
equation by multiplying with $x/\x$, that is, with the quotient of the
entries in a suitable column of the matrix defining the scroll.  We make
here a specific choice; we recall the notation $\lceil r \rceil$ for the
round up and $\lfloor r \rfloor$ for rounding off a real number $r$, and
note that $\auf m =\ab{m+1}$ for integers $m$.   This gives the following
equations:
$$
\displaylines{
\Rank \pmatrix{
z_0&z_1&\ldots&z_{k-1}&w_0&\ldots&w_{k-g-2}\cr
z_1&z_2&\ldots&z_k&w_1&\ldots&w_{k-g-1}\cr} \leq1,\cr
\phi_m=w_{\ab m}w_{\auf m}-\sumi a_iz_{\ab {m+i}}
z_{\auf {m+i}}, \qquad \hbox{$m=\mdot0{2k-2g-2}$.}
\cr}
$$

If $L$ is an arbitrary bundle of degree $d\geq2g+3$, then we can write
$L=kg^1_2+D$ with $\deg D= g+1-e$, and $k$ maximal such that  $D$ is
effective. The divisor $D$ can be described  by two polynomials $U(x)$
and $V(x)$, together with $F(x)=\sumi a_ix^i$ \cite{\rfTata, \S~1}.
Suppose $D=P_1+\cdots+P_{g+1-e}$ with all $P_i$ distinct and let
$(x_i,y_i)$ be the coordinates of $P_i$. Define $U(x)=\prod_i(x-x_i)$,
and take $V(x)$ the unique polynomial of degree $\le g-e$ with
$V(x_i)=y_i$. The ideal $(U,y-V)$ defines $D$; this is indeed a
subvariety of $y^2-F$, and one has $F-V^2=UW$ for some polynomial $W$.

The function $(y+V)/U$ defines a section of $H^0(C,L)$; a basis of this
vector space can be represented in inhomogeneous coordinates by the
polynomial forms:
$$
z_i=Ux^i, \qquad \hbox{$i=\mdot0k$;}\qquad
w_i=(y+V)x^i, \qquad \hbox{$i=\mdot0{k-e}$.}
$$
Modulo the equation $y^2-F$ we have the relation $(y+V)^2= UW+2(y+V)V$.
{}From it we obtain  $d-2g-1$ further equations $\phi_m$,  by rolling
factors. Write $U(x)=\sumj{g+1-e} U_ix^i$, $V(x)=\sumj{g-e} V_ix^i$ and
$W(x)=\sumj{g+1+e} W_ix^i$; this defines $F$. Then we have the following
equations:
$$
\displaylines{
\Rank \pmatrix{
z_0&z_1&\ldots&z_{k-1}&w_0&\ldots&w_{k-e-1}\cr
z_1&z_2&\ldots&z_k&w_1&\ldots&w_{k-e}\cr} \leq1,\cr
\phi_m=\sumj{g+1-e} U_iw_{\ab {m+i}}w_{\auf {m+i}}
-2\sumj{g-e} V_iw_{\ab {m+i}} z_{\auf {m+i}}
-\sumj{g+1+e} W_iz_{\ab {m+i}} z_{\auf {m+i}},\hfill\cr
\hfill \hbox{$m=\mdot0{d-2g-2}$.\qquad}
\cr}
$$

\roep (2.3)Rolling factors format {\rm \cite{\rfr}}.
This format occurs often in connection with divisors on scrolls. We start
with a  $k$-dimensional rational normal scroll $S\in\P=\P^n$; the
classical construction is to take $k$ complementary linear subspaces
$L_i$, spanning $\P$, with a parametrised rational normal curve
$\phi_i\colon\P^1 \to C_i\subset L_i$ of degree $d_i=\dim L_i$ in it, and
to take for each $p\in\P^1$ the span of the points $\phi_i(p)$. The
degree of $S$ is $f=\sum d_i$. The scroll $S$ can be defined by the
minors of the matrix:
$$
\def\zz#1#2#3{z^{(#1)}_{#2}&\ldots&z^{(#1)}_{#3}}
\pmatrix{
\zz10{d_1-1}&\ldots&\zz k0{d_k-1}\cr
\zz11{d_1}&\ldots&\zz k1{d_k}\cr}.
$$
The Picard group of $S$ is generated by the hyperplane class $H$ and the
ruling $R$ \cite{\rfschr, Sect.~1}; let $C$  be a divisor of type
$aH-bR$. The resolution of $\sier C$ as $\sier S$-module is $0\to \sier
S(-aH+bR) \to \sier S\to\sier C\to 0$. Schreyer describes, following
Eisenbud,  Eagon-Northcott type complexes ${\cal C}^b$ such that ${\cal
C}^b(a)$ is the minimal resolution of $\sier S(-aH+bR)$ as
$\sier\P$-module, if $b\geq-1$ \cite{loc.~cit., 1.2}. The resolution of
$\sier C$ is then obtained by taking a mapping cone. We obtain the
following first terms of the resolution:
$$
\displaylines{
\quad(\sier \P^f(-1)\otimes\Sym_{b-1}\sier \P^2 )(-a) \oplus
\wedge^3\sier \P^f(-1) \otimes \sier \P^2
\longrightarrow\hfill\cr\hfill
\longrightarrow(\Sym_b\sier \P^2 )(-a) \oplus
\wedge^2\sier \P^f(-1)
\longrightarrow \sier \P\longrightarrow \sier C\to 0.
\quad\cr}
$$
The  $b+1$ equations $\phi_m$, describing $C$ on $S$, are obtained by
`rolling' along $\P^1$: in terms of an inhomogeneous coordinate $x$ on
$\P^1$ we have $z^{(i)}_j=\psi^{(i)}x^j$ for some function $\psi^{(i)}$,
and in the transition from the equation $\phi_m$ to $\phi_{m+1}$, which
is obtained by multiplying with $x$, we have to increase the lower index
by one for exactly one of the factors $z^{(i)}_j$ in each monomial.

\roep (2.4) Infinitesimal deformations.
For the cone $X$ over a hyperelliptic curve of degree at least $2g+3$,
the dimension of $T^1_X(-1)$ is $2g+2$, by Proposition (1.25). We want to
describe the action on the equations. For this we have to study the
normal bundle of $C$.

Let  $d=\deg L\geq 2g+3$. For the scroll $S=\P(\pi_* L)$, we have $\sier
S(1)\sim E_0+af$, so the divisor of the line bundle $N_{C/S}(-1)$ is
$C\cdot (L+(2g+2-d)f)$. From the exact sequence:
$$
0\to H^0(N_{C/S}(-1))
\longrightarrow H^0(N_C(-1))
\longrightarrow H^0(N_S|_C(-1))
\longrightarrow H^1(N_{C/S}(-1))
$$
we obtain that $h^0(C,N_S|_C(-1))\leq h^0(C,N_C(-1))-
\chi(N_{C/S}(-1))=d+g+3-(3g+5-d)=2d-2g-2$, where $h^0(C,N_C(-1))$ can be
computed with the sequence  $ 0\longrightarrow  V^*\otimes H^0(C,\sier C)
\longrightarrow H^0(C,N_C(-1))\longrightarrow T^1_X(-1) \longrightarrow 0
$, by Theorem (1.4). On the other hand, $h^0(S,N_S(-1))=2d-2g-2$, and all
deformations are obtained by deforming the matrix of the scroll; because
$H^0(C,\sier C(1))=H^0(S,\sier S(1))$, the restriction map
$H^0(S,N_S(-1))\to H^0(C,N_S|_C(-1))$ is injective and so for dimension
reasons  an isomorphism. This shows that every infinitesimal deformation
of negative weight of $X$ comes from a deformation of the scroll.

We only compute the case $L=kg^1_2$; then $N_{C/S}(-1)=(2g+2-k)g^1_2$.

\beginprop(2.5)
Let $L=kg^1_2$ with $k>g+1$. On $T^1(-1)$ we  take
coordinates $(\mdot{s_1}{s_{k-1}},\mdot{t_0}{t_{2g+2-k}})$, subject to
the relations $\sumi a_is_{i+j}=0$ for $j=\mdot1{k-2g-3}$. The first
order deformations are the given by:
$$
\displaylines{
\Rank \pmatrix{
z_0&z_1+s_1&\ldots&z_{k-1}+s_{k-1}&w_0&\ldots&w_{k-g-2}\cr
z_1&z_2&\ldots&z_k&w_1&\ldots&w_{k-g-1}\cr} \leq1,\cr
\phi_m + \sum_{i=0}^{2g+2-k}t_iz_{m+i}
-\sum_{j=0}^{2g+2}\sum_{i=m+j}^{2g+2}a_is_{i-\auf j}z_{m+\auf j}
+\sum_{j=1}^{m-1}\sum_{i=0}^{m-j-1}a_is_{i+\ab
j}z_{m-\ab j},
\cr}
$$
where $\phi_m$ is as above, $m=\mdot0{2k-2g-2}$. Furthermore,  $a_j=0$
for $j>2g+2$, $s_0=0$ and $s_j=0$ for $j\geq k$. Therefore in the last
term the summation really runs to $\min(2g+2,m-j-1)$, and in the other
term to $\min(2g+2,k-1+\auf j)$.
\endprop

\bewijs
Consider the following relation, involving $\phi_m$:
$$
\displaylines{
\quad
z_{j+1}\phi_m-z_j\phi_{m+1}+w_\auf m(z_jw_\auf{m+1}
-z_{j+1}w_\ab m) \hfill\cr\hfill{}-\sumi
a_iz_\auf{m+i}(z_jz_\auf{m+i+1}-z_{j+1}z_\ab{m+i})=0.
\quad\cr}
$$
Now deform $\phi_m$ to $\phi_m+\phi_m'$, and the other equations as given
by the matrix. The condition that this defines an infinitesimal
deformation, is that in the local ring of the singularity the following
equation holds:
$$
z_{j+1}\phi_m'-z_j\phi_{m+1}'+\sumi
a_is_\ab{m+i}z_jz_\auf{m+i}=0.
$$
In terms of the inhomogeneous coordinates $(x,y)$ we can  write:
$$
x\phi_m'-\phi_{m+1}'+x\sumi
a_is_\ab{m+i}x^\auf{m+i}=0.
$$
This equation is independent of $j$. The relations involving
$w_{j+1}\phi_m$ give the same set of equations. It is easily checked that
the given formula for $\phi_m'$ satisfies this set of equations.
Furthermore, it gives a well defined element of the local ring: we check
that no $z_i$ with $i>k$ occurs with non-zero coefficient. If in the last
term $m-j-1<2g+2$, then $\ab j\geq \ab m -g -1$ and $m-\ab j\leq \auf m
+g+1\leq k$, because $m\leq 2k-2g-2$. Likewise, if $m+j\leq 2g+2$, then
$m+\auf j\leq k$.
\qed

\roep (2.6)$T^2$.
To determine the dimension of $T^2$ we use the {\sc Main Lemma} of
\cite{\rfbc, 1.3.2}. Let $h\colon X\to \C$ define a hyperplane section
$Y=h^{-1}(0)$ of $X$, then $\dim T^2_X/hT^2_X= \dim T^1_Y-e_h$, where
$e_h$ is the dimension of the smoothing component, on which the smoothing
$h$ of $Y$ lies. A general hyperplane section $Y$ is the cone over $d$
points in $\P^{d-g-1}$, lying on a rational normal curve of degree
$d-g-1$. Equations for this curve singularity can of course also be
written in rolling factors format: let the polynomial $F(\x,x)=\sumj d
a_i\x^{d-i}x^i$ determine the points on the rational curve, then we have:
$$
\displaylines{
\Rank \pmatrix{
z_0&z_1&\ldots&z_{d-g-2}\cr
z_1&z_2&\ldots&z_{d-g-1}\cr} \leq1,\cr
\phi_m=\sumj d a_iz_{\ab {m+i}}
z_{\auf {m+i}}, \qquad \hbox{$m=\mdot0{d-2g-2}$.}
\cr}
$$

\beginlem(2.7)
The $\sier Y$-module $T^2_Y$ is annihilated by the maximal ideal
${\frak m}_Y$.
\endlem

\bewijs
Let $N= {d-g-1\choose 2}+d-2g-1$ be the number of equations. We want to
show that $h\psi\in{\rm Im} \Hom((\sier X)^{\oplus k}, \sier X)$ for
every $\psi\in \Hom(R/R_0,\sier X)$ and every $h\in{\frak m}_X$. In the
determinantal case this is shown in \cite{\rfbc, 2.1.1}; we remark that
the proof as written is not correct: some equalities do not hold modulo
$R_0$, but only modulo the larger submodule $R_I$ of relations with
entries in the ideal $I$ of $X$; this is not a serious problem, because
$\Hom(R/R_I,\sier X)=\Hom(R/R_0,\sier X)$. Using their computations we
may assume that $h\psi$ vanishes on determinantal relations.

To describe the additional relations, we introduce the notation
$f_{i,j}=z_iz_{j+1}-z_{i+1}z_j$. We get:
$$
R_{j,m}=\phi_{m+1}z_j-\phi_m z_{j+1}-
\sumj d a_if_{j,\ab {m+i}}z_{\auf {m+i}},
$$
where $0\le j<d-g-1$ and $0\le m< d-2g-2$. The determinantal
relations are $R_{i,j,k}=f_{i,j}z_k-f_{i,k}z_j+f_{j,k}z_i$ and
$S_{i,j,k}=f_{i,j}z_{k+1}-f_{i,k}z_{j+1}+f_{j,k}z_{i+1}$.
We have the following equality:
$$
R_{j,m}z_k-R_{k,m}z_j-
\sumj d a_iR_{j,k,\ab {m+i}}z_{\auf {m+i}}=
\phi_mj(f_{j,k})-f_{j,k}j(\phi_m).
$$
Now fix an element $z_j\in{\frak m}_Y$. We look for $\psi_m\in \sier Y$
with $\psi(z_jR_{k,m})=\psi_{m+1}z_k-\psi_mz_{k+1}$ for all $k$ and $m$.
Because $\psi(z_jR_{k,m})=z_k\psi(R_{j,m})$, we can determine the
$\psi_m$ from the $\psi(R_{j,m})$, with the equations
$\x\psi(R_{j,m})-\x\psi_{m+1}+x\psi_m$; take $\psi_j=0$ and solve. For
$m>j$ we get that $\x^{m-j}\psi_{j,m} = \x^{m-j-1}x\psi(R_{j,m-1})+\cdots
+x^{m-j}\psi(R_{j,j})$; because  $z_{j+k}\psi(R_{j,l})=
z_j\psi(R_{j+k,l})$ we can divide by $\x^{m-j}$. A similar argument shows
that we can solve for $\psi_m$ with $m<j$.
\qed

Important for our application is that the corresponding result holds for
$T_X^2$; the computation is similar to the one for the hyperplane
section, but the formula's are more complicated. I expect that there is a
general statement for the rolling factors format, but it is not quite
clear what the `generic' rolling factors  singularity is.

We introduce more notation:
$g_{i,j}=z_iw_{j+1}-z_{i+1}w_j$
and $h_{i,j}=w_iw_{j+1}-w_{i+1}w_j$. We have a relation
$R_{j,m}$:
$$
\displaylines{
\quad R_{j,m}=\phi_{m+1}z_j-\phi_m z_{j+1}-
\sum U_ig_{j,\ab {m+i}}w_{\auf {m+i}}
+2\sum V_if_{j,{m+i\over2}} w_{m+i\over2}\hfill\cr
\hfill{}+2\sum V_ig_{j, {m+i-1\over2}} z_{m+i+1\over2}
+\sum W_if_{j,\ab {m+i}} z_{\auf {m+i}},\quad\cr}
$$
where one term with $V_i$ has to be chosen, depending on the parity of
$m+i$; we have similar relations $S_{j,m}$, involving $w_j$ and
$w_{j+1}$.  We find as above that $R_{j,m}z_k-R_{k,m}z_j\equiv 0
\pmod{R_D}$, where $R_D$ is the submodule, generated by the trivial and
determinantal relations. In the same way $S_{j,m}z_k-R_{k,m}w_j\equiv
0\pmod{R_D}$ and $S_{j,m}w_k-S_{k,m}w_j\equiv 0\pmod{R_D}$.  This gives
the desired result:

\beginlem (2.8)
The module $T^2_X$ is annihilated by the maximal ideal.
\endlem

\beginprop (2.9)
Let $Y$ be the cone over $d$ points on the rational normal curve of
degree $d-g-1$, with $d<2(d-g-1)$. Then $\dim\tee{-1}=d$,
$\dim\tee0=(g-1)(d-g-1)$, and $\dim\tee\nu=0$ for $\nu\neq0,-1$.
\endprop

\bewijs
The $d$ points are in generic position \cite{\rfgm, 3.3}: a set of $d$
points in $\P^k$ is in {\sl general position\/}, if for every $n$ the
images of the points under the $n$-tuple Veronese embedding of $\P^k$
span a linear space of maximal possible dimension, i.e. of dimension
$\min\{d,{k+n\over n}\}-1$. For points on a rational normal curve of
degree $k$ in $\P^k$ we have to consider the composed $kn$-tuple
embedding of $\P^1$; the dimension of the span is $\min\{d,kn\}-1$. In
our case $d<2k$, so every subset of the $d$ points is in general
position, and by definition the $d$ points are in {\sl generic
position\/}. Therefore $Y$ is negatively graded in the sense of
\cite{\rfpi},  cf. \cite{\rfgm, 3.3}, so $\dim\tee\nu=0$ for $\nu>0$.

The dimension of $\tee0$ is equal to the number of moduli for $d$ points
in $\P^{d-g-1}$. Because every (quadratic) equations occurs in a linear
relation,  $T^1_Y$ vanishes in degree $<-1$. We finish the proof by
computing $\tee{-1}$ in three steps.

\medbreak\noindent {\sl Step 1. \enspace
Every infinitesimal deformation of $Y$ of degree $-1$ comes from a
deformation of the cone over the rational normal curve of degree
$d-g-1$.\/}

\bewijs We may assume that $(\x,x)=(1,0)$ or $(0,1)$  is not a root of
the polynomial $F$, i.e. $a_0\ne0$ and $a_d\ne0$. Then $\x$ and $x$ are
not zero divisors in $\sier Y$. Therefore we get from the determinantal
relations the equations:
$$
f_{j,j+k}\x^{k-1}x=f_{j+1,j+k}\x^k+f_{j,j+1}x^k.
$$
By induction we find:
$$
f_{j,j+k}\x^{k-1}x^{k-1}=\sumj{k-1}
f_{j+i,j+i+1}\x^{2i}x^{2(k-i-1)}.   \eqno(*_{j,k})
$$
Let $n\in\Hom(I/I^2,\sier Y)$ be a normal vector of degree $-1$ with
$n(f_{i,i+1})=g_i$. From equation $(*_{j,k})$ we get that
$n(f_{j,j+k})\x^{k-1}x^{k-1}=\sum g_{j+i}\x^{2i}x^{2(k-i-1)}$. We claim
that $g_{j+i}$ is divisible by $\x^{k-i-1}x^i$. Every $g_{j+i}$ occurs in
$(*_{j+1,k})$ or $(*_{j,k-1})$, so by induction we obtain
$g_j=\x^{k-2}g_j'$, $g_{j+k-1}=x^{k-2}g_{j+k-1}'$ and the claim for
$0<i<k-2$. Therefore we have $x^2g_j'+\x^2g_{j+k-1}'-\x x\psi\equiv 0
\pmod F$ for some polynomial $\psi(\x,x)$.  The degree of $x^2g_j'(\x,x)$
is at most $d-g+1$; because $g>1$, we have in fact a polynomial equation
in $\x$ and $x$, so $g_j'$ is divisible by $\x$, and $g_{j+k-1}'$ by $x$.
By taking $k$ maximal we see that $g_j=\lambda_jz_j+\mu_jz_{j+1}
+\nu_jz_{j+2}$ for some $\lambda_j,\mu_j,\nu_j\in\C$: this is the formula
for the infinitesimal deformations of the cone over the rational normal
curve.

\medbreak\noindent {\sl Step 2. \enspace
The following formula's define a $d$-dimensional subspace of
$\tee{-1}$; deformation parameters are $(\mdot{s_1}{s_{d-g-2}},
\mdot{t_0}{t_{g+1}})$.
$$
\displaylines{
\Rank \pmatrix{
z_0&z_1+s_1&\ldots&z_{d-g-2}+s_{d-g-2}\cr
z_1&z_2&\ldots&z_{d-g-1}\cr} \leq1,\cr
\quad\phi_m + \sumj{g+1}t_iz_{m+i}+\sum_{i=1}^{d-g-2}
\sum_{j=0}^{2i-2m-1}
s_ia_{i+1+\auf j}z_{m+1+\auf j}
\hfill\cr
\hfill{}-\sum_{i=1}^{d-g-2}\sum_{j=0}^{2m-2i-1}
s_ia_{2i+1-m+\auf j}z_{i+1+\auf j},\quad
\cr}
$$
where $\phi_m$ is as above, $m=\mdot0{2k-2g-2}$.
Furthermore,  $a_j=0$ for $j>d$ or $j<0$.}

\bewijs
The deformation of the matrix gives the equation $f_{j,k}+
s_jz_{k+1}-s_kz_{j+1}$. Suppose $\phi_m$ is deformed to $\phi_m+\phi_m'$.
The relations $R_{j,m}$ give:
$$
\x\phi_{m+1}'-x\phi_m'+\sumj da_is_\ab{m+i}x^{\auf{m+i}+1}=0.
$$
The given formula satisfies these equations. The indices in
it do not exceed $d-g-1$.

\medbreak\noindent {\sl Step 3. \enspace
The deformations of the scroll, which deform $f_{j,j+1}$ by
$g_j=\mu_jz_{j+1}$, do not extend to deformations of $Y$.\/}

\bewijs
Suppose $f_{i,j}$ and $\phi_m$ deform with $g_{i,j}$ and $\phi_m'$. From
equation $(*_{j,k})$ we get that  $g_{j,j+k}=\sum_{i=j}^{j+k-1}
\mu_ix^{2j+k-i}$. The relation $R_{j,m}$ gives that
$$
\phi_{m+1}'-x\phi_m'+\sumj {2j-m-1}a_i \sum
_{k=\ab{m+i}}^{j-1} \mu_kx^{m+i-k}-
\sum_{i=2j+2-m}^da_i\sum_{k=j}^{\ab{m+i}-1} \mu_kx^{m+i-k}=0.
$$
Because $F(\x,x)=0$, this expression is independent of $j$,
so we take $j=0$. Let $\alpha$ be the smallest index such
that $\mu_\alpha\neq0$. We set $m=2d-2g-3$:
$$
\phi_{2d-2g-2}'-x\phi_{2d-2g-3}'-
\sumj d \sum_{k=g+3-\auf d}^{i-\alpha}
a_i\mu_{i-k}x^{2d-2g-3+k}=0.
$$
Because $d-g-1$ is the highest power of $x$, which can occur in
$\phi_{2d-2g-3}'$, the coefficient of $3d-2g-3-\alpha\geq 2d-g-1$ in the
sum vanishes; so $a_d\mu_\alpha=0$, which contradicts our assumption that
$a_d\neq0$. Therefore all $\mu_i$ are zero. \qed

\beginprop(2.10)
Let $X=X(C,L)$ be the cone over a hyperelliptic curve of genus $g$,
embedded with a complete linear system $L$ of degree $d$, with $d>2g+3$.
Then
$\dim\ttw{-2}=d-2g-3$,
$\dim\ttw{-1}=(g-2)(d-g-3)$, and
$\dim\ttw\nu=0$ for $\nu\neq -1,-2$.
\endprop

\bewijs
For the hyperplane section $Y$ the previous Proposition gives $\dim
T^1_Y=d+(g-1)(d-g-1)$; the dimension of a smoothing component is
$e=\mu+t-1$ \cite{\rfgm, 2.5.(3)}, where $t$ is the type of the
singularity. In our case $t=g=\delta(Y)-d+1$, so $e=3g+d-2$. Because
$T^2_X$ is annihilated by the maximal ideal, the {\sc Main Lemma} of
\cite{\rfbc} gives $\dim T^2_X=\dim T^1_Y-e=(g-1)(d-g-4)-1$. By the
formula of  (1.26) we have  $\ttw{-1}=(g-2)h^0(L)-6(g-1)+\dim T^1_X(-1)$.
By Lemma (1.28), $\dim\ttw\nu=0$ for $\nu< -1$.
\qed

\roep(2.11) Remark.
One can also compute $\ttw{-2}$ directly; we sketch this here for the
special case $L=kg^1_2$. We have to show that  $\dim\ker\{H^0(L)\otimes
H^1(L^{-1}) \to H^1(\sier C)^{\oplus N}\} = 10k+2g-7$, where $N$ is the
number of equations. The elements $y/(\x^{g+1-i}x^i)$, $i=\mdot 1g$, form
a basis of $H^1(\sier C)$; we abbreviate $[i]=y/(\x^{g+1-i}x^i)$, for all
$i\in {\Bbb Z}$. A basis of $H^0(L)\otimes H^1(L^{-1})$ is:
$$
\def\part#1#2#3{{#1\over
\x^{#2-i}x^i}{\partial\over\partial #3_j}}
\openup1\jot
\tabskip=0pt plus1fil
\halign to\displaywidth
{\tabskip=0pt
$\hfil#_{i,j}$&$\displaystyle{}=#\hfil
$&\qquad$i=\mdot1#,\hfil $&\quad$j=\mdot0#\hfil$
\tabskip=0pt plus1fil\cr
\phi&\part 1kz & k-1& k ,\cr
\psi&\part y{g+k+1}z& g+k& k, \cr
\chi&\part 1kw & k-1& k-g-1, \cr
\theta&\part y{g+k+1}w& g+k& k-g-1.
\cr }
$$
On $f_{l,m}$ only the action of $\psi_{i,j}$ is non trivial:
\def\act#1={\psi_{i,#1}(f_{l,m})=}
$\act l=[i-m-1]$, $\act l+1=-[i-m]$, $\act m+1=[i-l]$, $\act m=[i-l-1]$
and $\act j=0$ otherwise.  Therefore:
$$
\sum_{i,j}\alpha_{i,j}\act j=
\sum_i (\alpha_{i+m+1,l}+\alpha _{i+l,m+1}
-\alpha_{i+m,l+1}-\alpha_{i+l+1,m})[i].
$$
The coefficient of $[i]$ has to vanish for $1\leq i \leq g$. A set of
solutions to these equations is given by $\alpha_{i,j}=\alpha_{i+j}$ for
some constants $\alpha_j$, $j=\mdot1{g+2k}$. Define
$\beta_{i,j}=\alpha_{i+1,j}-\alpha_{i,j+1}$, where $i=\mdot1{g+k-1}$ and
$j=\mdot0{k-1}$. We have the equations $\beta_{i+m,l}=\beta_{i+l,m}$ for
$i=\mdot1g$ and $0\leq l<m\leq k-1$. Consider pairs $(i,j)$ and $(i',j')$
with $i+j=i'+j'$, and suppose that $j=j'+s(g-1)+(g-r)$. Write
$(i,j)\sim(i',j')$ if $\beta_{i,j}=\beta_{i',j'}$. If $i-j'>g$, we can
use the equations to obtain $(i,j)\sim (t+j,i-t)\sim(i+g-t,j-g+t)$, where
$t=r$ if $s=0$, or $t=1$ for $s>0$; for $i-j'<0$ we have a similar
inductive procedure, whereas otherwise a suitable equation directly gives
that $(i,j)\sim(i',j')$ if  $i+j=i'+j'$.

For the action of $\chi_{i,j}$ on $h_{l,m}$ we have analogous
computations. The  next step is to compute the action on $\phi_m$:
$$
\def\Auf#1{\mathchoice{{\textstyle\auf{#1}}}%
{{\scriptstyle\auf{#1}}}{{\scriptscriptstyle\auf{#1}}}{}}
\def\Ab#1{\mathchoice{{\textstyle\ab{#1}}}%
{{\scriptstyle\ab{#1}}}{{\scriptscriptstyle\ab{#1}}}{}}
\displaylines{\quad
\sum ((\alpha_{i+j}+j\beta_{i+j})\psi_{i,j}
+(\gamma_{i+j}+j\delta_{i+j})\chi_{i,j})(\phi_m)=
\hfill\cr\qquad{}=
\sum(\gamma_{i+\Ab m}+\Ab m\delta_{i+\Ab m})[\Auf m-i]
+(\gamma_{i+\Auf m}+\Auf m\delta_{i+\Auf m})[\Ab m-i]
\hfill\cr
\qquad\qquad{}-a_l (\alpha_{i+\Ab {m+l}}+\Ab {m+l}
\beta_{i+\Ab {m+l}})[\Auf {m+l}-i]
\hfill\cr
\qquad\qquad{}-a_l(\alpha_{i+\Auf {m+l}}+\Auf {m+l}
\beta_{i+\Auf {m+l}})[\Ab {m+l}-i]\hfill\cr
\qquad{}=
\sum(2\gamma_{m-c}+m\delta_{m-c}-2a_l\alpha_{m+l-c}
-a_l(m+l)\beta_{m+l-c})[c]\hfill\cr}
$$
We get equations $\delta_j-\sumi a_i\beta_{j+i}$ and $2\gamma_j-\sumi
a_i(2\alpha_{j+i}-i\beta_{j+i})$, except for $j=2k-g-2$, when there is
only one equation; but there is no $\delta_{2k-g-2}$. So the parameters
$\gamma_j$ and $\delta_j$ are completely determined by  $\alpha_j$ and
$\beta_j$. As $\beta_1=\beta_{2k+g}=0$, we have a $(4k+2g-2)$-dimensional
solution space.

Finally, we compute the action on $z_lw_{m+1}-z_{l+1}w_m$:
$$
\sum(\alpha_{i,j}\phi_{i,j}+\beta_{i,j}\theta_{i,j})(g_{l,m})=
\sum_i (\alpha_{i+m+1,l}+\beta _{i+l,m+1}
-\alpha_{i+m,l+1}-\beta_{i+l+1,m})[i].
$$
The coefficient of $[i]$ certainly vanishes, if $\alpha_{i,j}$
and $\beta_{i,j}$ depend only on $i+j$.  Define
$\gamma_{i,j}=\alpha_{i+1,j}-\alpha_{i,j+1}$ and
$\delta_{i,j}=\beta_{i+1,j}-\beta_{i,j+1}$. Then
$\gamma_{i+m,l}=\delta_{i+l,m}$. One shows that
$\gamma_{i,j}=\delta_{i,j}=\epsilon_{i+j}$ for all $i$, $j$.
The solution space has dimension $6k-5$. \qed

The dimension of $\ttw{-2}$ is two less than the number of equations
$\phi_m$. This fact has an interpretation in terms of the  rolling
factors format, if all the $\phi_m$ are quadratic equations.

\beginprop(2.12)
Let $X$  be the cone over a divisor of  type  $2H-bR$ on a scroll $S$.
Consider infinitesimal deformations of negative degree of $X$, which come
from  deformations of the cone over $S$; if the versal deformation of the
cone over $S$ is:
$$
\def\zz#1#2#3{z^{(#1)}_{#2}&\ldots&z^{(#1)}_{#3}}
\pmatrix{
z_0^{(1)}&z^{(1)}_{1}+s^{(1)}_{1}&\ldots
&z^{(1)}_{d_1-1}+s^{(1)}_{d_1-1}&\ldots&\zz
k0{d_k-1}+s^{(k)}_{d_k-1}\cr
z^{(1)}_1&\zz12{d_1}&\ldots&\zz k1{d_k}\cr},
$$
then there are additional parameters $t_i$, and linear equations on the
$t_i$ and $s_i^{(j)}$, such that the additional equations are given by
$\phi_m+\phi_m'( \boldmath t,\boldmath s, \boldmath z)$, $m=\mdot0b$.
Then this deformation can be extended over a base space given by $b-1$
equations $\phi_m'(\boldmath t,\boldmath s, \boldmath s)-\phi_m(\boldmath
s)=0$, $m=\mdot1{b-1}$.
\endprop

\bewijs
We have to lift the relations involving the $\phi_m$. The rolling factors
assumption gives that we can write
$\phi_m=\sum_{\alpha}z_{\alpha}c_{\alpha}$ and
$\phi_{m+1}=\sum_{\alpha}z_{\alpha+1}c_{\alpha}$, where $c_\alpha$
depends (linearly) on $\boldmath{z}$, and $\alpha$ runs through all
possible indices ${}^{(i)}_j$, and $c_\alpha$ depends (linearly) on
$\boldmath{z}$. This gives the relation
$\phi_{m+1}z^{(i)}_j-\phi_mz^{(i)}_{j+1}- \sum
f_{j\alpha}^{(i)}c_\alpha$, where $f_{j\alpha}^{(i)}$ is the
determinantal equation $z^{(i)}_jz_{\alpha+1}-z^{(i)}_{j+1}z_\alpha$. We
lift it  to:
$$
\displaylines{\quad
(\phi_{m+1} +\phi_{m+1}'(\boldmath z))(z^{(i)}_j+s^{(i)}_j) -(\phi_m
+\phi_m'(\boldmath z+\boldmath s))z^{(i)}_{j+1} -\sum
\tf^{(i)}_{j,\alpha}c_\alpha\hfill\cr \hfill\equiv \sum z^{(i)}_{j+1}
s_\alpha c_\alpha(\boldmath s) \pmod {\cal I},\quad\cr}
$$
where ${\cal I}$ is the ideal of the deformed scroll, and
$\tf_{j\alpha}^{(i)}$ is a deformed equation. The lift up to first order
is possible by assumption. If $1\leq m\leq b-1$, then $\phi_m$ occurs in
a relation as first and as second term. Therefore  $\phi_m'(\boldmath z)$
and $\phi_m'(\boldmath z+\boldmath s)-\phi_m(\boldmath z)$ have to be
equal.
\qed

\roep (2.13)Example.
Let $X$ be the cone over a hyperelliptic curve, embedded with
$(2g+2)g^1_2$. With Proposition (2.5) we find as  equations
for the negative degree part of the versal base:
$$
ts_m-
\sum_{\scriptstyle i=m+j \atop\scriptstyle j=0}
 ^{\scriptstyle j=2g+2\atop\scriptstyle i=2g+2}
a_is_{i-\auf j}s_{m+\auf j} +
\sum_{\scriptstyle i=0\atop\scriptstyle j=1}
 ^{\scriptstyle j=m-1\atop\scriptstyle i=m-j-1}
a_is_{i+\ab j}s_{m-\ab j}
+\sumi a_is_{\auf{m+i}}s_{\ab{m+i}}.
$$
In particular, if we take  the curve $y^2-1+x^{2g+2}$, we have:
$$
ts_m+\sum_{j=1}^{2g+1-m}s_{2g+2-\auf j}s_{m+\auf j}
+\sum_{j=2}^{m}s_{\ab j}s_{m-\ab j},
\qquad
m=1\hbox{, \dots, }2g+1.
$$
For
$g=2$ there are  five equations:
$ts_1+2s_5s_2+2s_4s_3$,
$ts_2+2s_5s_3+s_4^2+s_1^2$,
$ts_3+2s_5s_4 +2s_1s_2 $,
$ts_4+s_5^2 +2s_1s_3+s_2^2$,
$ts_5+2s_1s_4+2s_2s_3 $.
This is a complete intersection of degree $2^5$: it is the cone
over 32 distinct points.

We shall prove in the next Section that the base space is always a
complete intersection. It is difficult to see this from the above
equations.


\beginsection 3. Smoothing components

(3.1) In this section we prove that the cone over a hyperelliptic curve of
degree $4g+4$ has $2^{2g+1}$ smoothing components. To each component
corresponds a surface with the curve as hyperplane section; these are
ruled surfaces, and we show how to obtain all by elementary
transformations on a given one.  To this end we identify the one
dimensional subspace in $T^1$, determined by the surface.

Let $\wlx\subset\P^{N+1}$ be the projective cone over a projectively
normal algebraic variety $C\subset\P^N$ of dimension $n$, and let
$S\subset\P^{N+1}$ be variety of dimension $n+1$ with $C$ as hyperplane
section. Then there exists a one parameter deformation of $\wlx$ with $S$
as general fibre: let $\wly\subset\P^{N+2}$ be the projective cone over
$S$, and consider the pencil $\{H_t\}$ of hyperplanes through $\P^N
\supset C$; let  the hyperplane $H_0$ pass through the vertex, then
$\wly\cap H_0=\wlx$, while for $t\neq0$ the projection from the vertex
establishes an isomorphism between $\wly\cap H_t$ and $S$. Pinkham calls
this construction `sweeping out the cone' \cite {\rfpi}; the degeneration
to the cone occurs already in the famous `Anhang F' \cite{\rfse}.

For the affine cone $X$ over $C$ we have a deformation  with Milnor fibre
$S-C$. The versal base $S_X^-$ in negative degree is a fine moduli space
for so called $R$-polarised schemes \cite{\rflo}, a notion defined in
general for quasi-homo\-ge\-ne\-ous spaces; in our situation this are
spaces $S$ with $C$ as hyperplane section. Basically one considers the
coordinate $t$, which defines the hyperplane section, as deformation
parameter.

Let $S$ be given by equations $F_i({\boldmath x},t)=f_i({\boldmath
x})+tf_i^{(1)}+\dots+t^{d_i}f_i^{(d_i)}$, in homogeneous coordinates
$({\boldmath x},t)$, where $\deg f_i^{(j)}=d_i-j$.  The base space of $X$
has a $\C^*$-action. For simplicity we assume that the only occurring
negative  degree is $-1$, as in our application. The by $S$ induced
infinitesimal deformation of $X$ is given by $f_i\mapsto {\partial\over
\partial t}F_i({\boldmath x},t)|_{t=0}=f_i^{(1)}$. This can be
interpreted as section of $H^0(C,N_C(-1))$.

Let $S$ be embedded by the line bundle $L$, and let $\P=\P(H^0 (S,L)^*)$.
Let $C=H\cap S$, with $H$ a hyperplane in $\P$, and let $i\colon C \to S$
be the inclusion. Then $N_{C/S}=L_{|_C}$. In the exact sequence
$$
0\longrightarrow \sier C\cong
N_{C/S}(-1)\longrightarrow N_{C/\P}(-1) \longrightarrow
i^*N_{S/\P}(-1) \longrightarrow 0
$$
the infinitesimal deformation is the image of $1\in\sier  C$. The curve
$C$ is minimally embedded in the hyperplane $H$, so we really want a
section of $N_{C/H}(-1)$; the exact sequence
$$
0\longrightarrow
N_{C/H}(-1)\longrightarrow N_{C/\P}(-1) \longrightarrow
i^*N_{H/\P}(-1) \longrightarrow 0
$$
splits: $N_{H/\P}(-1)$ is
generated by the global section, which sends the equation $t$ to 1,
and this section can be mapped to the section $F\mapsto
{\partial\over \partial t}F|_{t=0}$ of $N_{C/\P}(-1)$.

For  curves we have a third description, which uses $(T^1)^*$. We recall
from Proposition (1.3) the  sequence $0\longrightarrow N^*_C\otimes L
\longrightarrow  V\otimes_\C \sier C  \longrightarrow {\cal
P}^1_C(L)\longrightarrow 0$. We have seen that $T^1_X(-1)^*=\coker
\{V\otimes H^0(K) \to H^0({\cal P}^1_C(L)\otimes K)\}$.  Consider the
following diagram of  exact sequences:
$$
\cdmatrix{
&&&&0&&0&&\cr
&&&&\Big\downarrow&&\Big\downarrow&&\cr
0 & \longrightarrow & \sier C &
\longrightarrow & \Omega^1_S\otimes \sier C(1)
& \longrightarrow & \Omega^1_C (1) & \longrightarrow & 0 \cr
 & & \Big\Vert & &\Big\downarrow & & \Big\downarrow & & \cr
 0 & \longrightarrow & \sier C & \longrightarrow &
{\cal P}_S^1(L)\otimes\sier C &
 \longrightarrow & {\cal P}^1_C(L) & \longrightarrow & 0 \cr
&&&&\Big\downarrow&&\Big\downarrow&&\cr
&&&&L&\Relbar\joinrel\Relbar&L&&\cr
&&&&\Big\downarrow&&\Big\downarrow&&\cr
&&&&0&&0&&\cr
}
$$

\beginprop(3.2)
Let $\xi\in H^0(C,N_C(-1))$ be the by $S$ induced infinitesimal
deformation of the cone over $C$. The  map $\cup\xi\colon T^1(-1)^*\to
\C$ is  the connecting homomorphism $H^0({\cal P}^1_C(L)\otimes K)\to
H^1(K)$, obtained from the above sequence by tensoring with
$K=\Omega^1_C$; alternatively one may consider the map $H^0(K^2L)\to
H^1(K)$.
\endprop

\bewijs
We first describe $H^0({\cal P}^1(L)\otimes K)$. The sheaf ${\cal
P}^1(L)$ is generated by global sections; if we take a local coordinate
$u$ on $C$ and write a Newton dot for the derivative with respect to $u$,
then the map $V\otimes\sier C\to{\cal P}^1(L)$ is given in accordance
with our earlier notations by $dz_i\mapsto z_i(u)+\dot z_i(u)du$. Let
$\phi$ be a global section of  ${\cal P}^1(L)\otimes K$. On $U_j=\{
z_j\neq0\}$ the section $\phi$ can be represented by $\sum
\omega_j^idz_i$, and  the representations for different $j$ are connected
by the conditions  $\sum \omega_j^iz_i=\sum \omega_k^iz_i$ and  $\sum
\omega_j^i\dot z_i=\sum \omega_k^i\dot z_i$.

This allows us to compute $\phi\cup \xi$: the cochain  $\sum
(\omega_k^i-\omega_j^i)dz_i$ represents an element of
$H^1(N^*_C(1)\otimes \Omega_C^1)=H^0(N_C(-1))^*$. To express $\xi$ in
terms of the $\partial\over\partial z_i$, we take local coordinates
$(u,\tau)$ on $S$ with $\tau$ vanishing on $C$, and we denote
differentation w.r.t. $\tau$ by $'$. For every equation $F$ of $S$ we
have $\partial F/\partial \tau=0$, so:
$$
t'(u,0){\partial F\over \partial t}+\sum z_i'(u,0){\partial F\over
\partial z_i}=0,
$$
and  therefore $\xi=-\sum (z_i'/t'){\partial\over \partial z_i}$.
Now $\phi\cup \xi=\sum(\omega_j^i-\omega_k^i)(z_i'/t')$.

To compute the connecting homomorphism $\delta$ in the exact
sequence:
$$
0\to \Omega^1_C\to {\cal P}^1_S(L)\otimes \Omega^1_C
\to {\cal P}^1_C(L)\otimes \Omega^1_C\to 0,
$$
we lift $\phi$ on $U_j$ to ${\cal P}^1(L)\otimes K$. On $U_j=\{
z_j\neq0\}$ we write $\phi$ as $\sum \omega_j^idz_i+\omega_j dt$, and the
cochain conditions  on $S$ are: $\sum \omega_j^iz_i+\omega_jt=\sum
\omega_k^iz_i+\omega_kt$, $\sum \omega_j^i\dot z_i+\omega_j\dot t=\sum
\omega_k^i\dot z_i +\omega_k\dot t$ and  $\sum
\omega_j^iz_i'+\omega_jt'=\sum \omega_k^iz_i'+\omega_kt'$; for $\tau=0$
we have $t=\dot t=0$, so the first two conditions are the same as before,
and we have $(\omega_k-\omega_j)t'=\sum(\omega_j^i-\omega_k^i)z_i'$. The
cochain $\omega_k-\omega_j$ represents $\delta(\phi)\in H^1(K)$, and
therefore $\delta(\phi)=\phi\cup\xi$.
\qed

We now concentrate on the hyperelliptic case. The cone over a
hypereliptic curve of degree $d\leq 4g+4$ is smoothable. More precisely,
we have:

\beginprop {\rm \ \cite{\rftenh}}(3.3)
Let  ${}^gC$ be a hyperelliptic curve, embedded with a complete linear
system $L$  of degree $7g/3+1\leq d\leq 4g+4 $. Then $C$ is a hyperplane
section of a projectively normal surface $\wl S=\phi_{|C|}(S)$, where
$S$  can be obtained from a rational ruled surface by blowing up $4g+4-d$
points. The dimension of the smoothing component is
$7g+4-d+h^0(C,K^2L^{-1})$. If   $d> \max\{4g-4,3g+6\}$ for $g\neq 6$, and
$d>25$ for $g=6$, then  every smoothing component is of this form, and
has dimension $7g+4-d$.
\endprop

In particular, for $d=4g+4$ the surface $S$ is ruled, and the dimension
of the smoothing component is $3g=(2g-1)+g+1$. The dimension of the
hyperelliptic locus is $2g-1$, and $g$ is the dimension of $\Pic^d$. So
for general $(C,L)$ each smoothing component determines a unique surface
$\wl S$ with $C$ as hyperplane section.  Furthermore $\dim
T^1_X(-1)=2g+2$, and the base space in negative degree is given by $2g+1$
equations.  From these facts we cannot yet conclude that the equations
define a complete intersection.

The existence of smoothings can be shown in the following way: given a
line bundle $L$ of degree $4g+4$ on $C$ (which is not necessarily the
linear system for which we want a smoothing --- we denote that
temporarily by $N$), we consider the  scroll $\wl S$ of type $(a,b)$, on
which $\phi_L(C)$ lies: it is the image of the ruled  surface
$S=\P(\pi_*L)$, cf. (2.2). On $S$ the hyperplane class is $L=E_0+af$, and
the curve $C$ is a divisor of type $2E_0+(g+1+e)f=2L-(2g+2)f$; here
$e=a-b$. So for every $L$ with $L^2\cong N\otimes (g^1_2)^{2g+2}$ the
normal bundle $N_{C/S}$ is isomorphic to  our fixed bundle $N$. The
number of solutions to this equation is the order of the group $J_2(C)$ of
2-torsion points on ${\rm Jac}(C)$,  which is $2^{2g}$.

For all surfaces $S$ obtained by this construction we have $e\equiv g+1
\pmod 2$, because $e=a-b$, and $a+b=3g+3$. This is not surprising,
because for a fixed curve and variable $L$ the surfaces $\P(\pi_*L)$ form
a continuous family, and the parity of $e$ is conservated under
deformations.

Let $C$ lie on on the ruled surface $S\cong F_e$, with normal bundle
$N_{C/S}=N$, with $0\leq e\leq g+1$.  Denote the elementary
transformation \cite{\rfHAG, V.5.1.7} of $S$ in $Q\in S$ by $\elm_Q$.
Then we have the following simple, but important observation:

\beginlem(3.4)
Let $C'$ be the strict transform of $C$ on $S'=\elm_P(S)$, where $P$
is a Weierstra\ss\ point. Then the normal bundle $N_{C'/S'}$ is
equal to the normal bundle $N_{C/S}$.
\endlem

This lemma shows the existence of ruled surfaces $C\subset S\cong
F_e$ with $N_{C/S}=N$ and $g+1+e$ odd; again  $C$ is of type
$2E_0+(g+1+e)f$.

\beginprop(3.5)
Let $B$ be the set of  Weierstra\ss\ points on $C$. Let $\elm_B$ be the
composition of the elementary transformations in all $P\in B$. Then
$\elm_B(S)$ is isomorphic to $S$, under an isomorphism $I$, which leaves
$C$ pointwise fixed; on the general fibre $f$ of $S\to\P^1$, $I$
restricts to the unique involution with $C\cap f$ as fixed points. The
isomorphims $I$ maps the linear system $|C|$ on $S$ to the linear system
$|C'|$ on $S'$.
\endprop

\bewijs
The rational map $\elm_B$ can be factorised as
$S\smash{\mathop{\longleftarrow}\limits^{\sigma}} \wt S
\mapright{\sigma'}S'$, where the maps $\sigma$ and  $\sigma'$ are the
blowups in the points of $B$, with exceptional curves $E_i$ and $E_i'$.
The involution on the general fibre f of $\wt S\to\P^1$ with $f\cap C$ as
fixed points extends to an involution $\wt I$ on $\wt S$, which
interchanges $E_i$ and $E_i'$.  This map descends to the required
isomorphism. The inverse image on $\wt S$ of the linear system $|C|$  is
$|\wl C+\sum E_i|$, where $\wl C$ is the strict transform of $C$. The
involution $\wt I$ transforms this system into $|\wl C+\sum
E_i'|=(\sigma')^*|C'|$.
\qed

\beginprop(3.6)
Let $C\subset S\cong F_e$ be a curve of type $2E_0+(g+1+e)f$ with normal
bundle $N_{C/S}=N$. Denote for a subset $T\subset B$ by $\elm_T$ the
composition of the $\elm_{P_i}$, $P_i\in T$. Two surfaces $\elm_{T_1}(S)$
and $\elm_{T_2}(S)$ induce the same deformation of the cone $X(C,N)$ if
and only if $T_1=T_2$ or $T_1=B\setminus T_2$.
\endprop

\bewijs
We have to compute the connecting homomorphism from the exact sequence:
$$
0\to \Omega^1_C\to \Omega^1_S\otimes \Omega^1_C(1)\to
(\Omega_C^1)^{\otimes 2}(1) \to 0,
$$
or, what amounts to the same, the hyperplane in $H^0(C,K^2(1))$, which is
the image of  $H^0(C,\Omega^1_S\otimes \Omega^1_C(1))$. We have on $S$
the sequence:
$$
0\to\pi^*\Omega^1_{\P^1}\to \Omega^1_S \to
\Omega^1_{S/\P^1}\to 0,
$$
and on $C$:
$$
0\to\pi^*\Omega^1_{\P^1}\to \Omega^1_C \to \Omega^1_{C/\P^1}\to
0,
$$
where $\Omega^1_{C/\P^1}=T^1_X(-1)^*\cong \oplus_{P\in B} \C_P$. The map
$\Omega^1_S\otimes \Omega^1_C(1)\to (\Omega_C^1)^{\otimes 2}(1)$ is an
isomorphism on the subspace $\pi^*\Omega^1_{\P^1} \otimes \Omega^1_C(1)$.
Because $K_S=-2E_0-(e+2)f$, we have $\Omega^1_{S/\P^1}=\sier S(-2E_0-ef)$,
and $\Omega^1_{S/\P^1}(C+(g-1)f)=\sier S(2gf)$. We identify $T^1_X(-1)^*$
with $H^0(C, (\Omega^1_C/\pi^*\Omega_{\P^1}) \otimes \Omega_C^1(1))$ by
taking in each Weierstra\ss\ point $P_i$ a  generator of $K^2(1)/{\frak
m}_{P_i}K^2(1)$; the map from $H^0(C,\Omega^1_{S/\P^1} \otimes
\Omega^1_C(1))$ consists in taking the coefficients in the points $P_i$.
The sections of $(g^1_2)^{\otimes2g}$ are $\mdot1{x^{2g}}$, which come
from sections on $S$, and $\mdot y{yx^{g-1}}$, which vanish in the
Weierstra\ss\ points. Therefore, if we take coordinates $s_i$ on
$T^1_X(-1)^*$, and if $x$-coordinate of $P_i$ is $a_i$, then the image
of  $H^0(C,\Omega^1_{S/\P^1} \otimes \Omega^1_C(1))$  is given by the
determinant:
$$
D(s_1,\dots,s_{2g+2})=
\left|\matrix{
1&a_1&a_1^2&\ldots&a_1^{2g}&s_1\cr
\vdots&\vdots&\vdots&\ddots&\vdots&\vdots\cr
1&a_{2g+2}&a_{2g+2}^2&\ldots&a_{2g+2}^{2g}&s_{2g+2}\cr
}\right|.
$$
Now we consider the surface $S'=\elm_T(S)$ for some $T\subset B$. Again
we can identify the sections of  $H^0(C,\Omega^1_{S'/\P^1}
\otimes\Omega^1_C(1))$, coming from $S'$, with the polynomials
$\mdot1{x^{2g}}$, but we have to express these in the same basis of
$T^1_X(-1)^*$, which we used for $S$. Consider local coordinates $(x,y)$
in a neighbourhood of a point $P\in T$, such that ruling is given by
$\pi(x,y)=x$, and $C$ is $y^2=x$. We blow $S$ up in $P$; the strict
transform of $C$ passes through the origin of the $(\eta,y)$ coordinate
patch, where $(x,y)=(\eta y,y)$. Now blow down the $y$-axis: we have
coordinates $(\xi,\eta)=(\eta y, \eta)$, so $x=\xi$, $y=\xi/\eta$.
Therefore the local generator $dy$ of $\Omega^1_{C,P}$ is transformed into
$d\eta$: we  have $dy=d\xi/\eta-\xi d\eta/\eta^2$; however, considered as
section of $\Omega^1_S/\pi^*\Omega^1_{\P^1}$, the formula makes sense for
$\eta\neq 0$ and is on $C$ the same as $-\xi d\eta/\eta^2=-d\eta$. This
computation shows that  $S'$ yields the hyperplane:
$$
D((-1)^{\chi_T(P_1)}s_1,\, \dots,\,
(-1)^{\chi_T(P_{2g+2})}s_{2g+2})=0,
$$
where $\chi_T$ is the
characteristic function of $T$, i.e~$\chi_T(P_i)=1$ if and only if
$P_i\in T$.

To finish the proof we remark that the coefficient of $s_i$ in the
equation $D$ is a Vandermonde determinant, and therefore non zero.
\qed

\beginthm(3.7)
Let $X$ be the cone over a hyperelliptic curve $C$ of genus $g$,
embedded with a complete linear system $L$ of degree $4g+4$.
Suppose $L\neq 4K$, if $g=3$. Then $X$ has $2^{2g+1}$ smoothing
components. \endthm

\bewijs
The number of subsets $T$ of $B$ modulo the equivalence relation $T\sim
B\setminus T$ is $2^{2g+1}$. So  Proposition (3.6) gives this number of
one parameter smoothings of $X$. After a rescaling these define in
$\P(T^1_X(-1))$  the points $(\pm 1 :\dots :\pm1)$, and this is a
complete intersection. In coordinates $s_i$ on $\P(T^1_X(-1))$  the ideal
of these points is generated by $s_i^2-s_j^2$. The $2g+1$ quadratic
equations for the base in negative degree are contained in this ideal;
they generate this ideal if and only if they are linearly independent.
Because the dimension of each smoothing component is $3g$, the equations
are independent for generic $(C,L)$. Suppose that for some special
$(C,L)$ the equations are dependent. Then the base space of $X(C,L)$ is
not smooth along at least one of the parameter lines  of the $2^{2g+1}$
smoothing we just constructed. But the only singularities in the fibres
are cones over a rational normal curve of degree $g+1$, which have a
smooth reduced base space, if $g\neq 3$. Therefore the equations are
always independent.

For $g=3$ we still have the same description of $S^-$, but in case $L=4K$
the cone over the rational normal curve of degree 4 appears as
singularity over one component. Its Veronese smoothing leads to an
additional smoothing component; it exists also for the non hyperelliptic
curves.
\qed

\roep (3.8) Remark.
The set $\cal B$ of subsets  $T$ of $B$ modulo the equivalence relation
$T\sim B\setminus T$ forms a group, isomorphic to ${\Bbb Z}^{2g+1}$; the
subgroup ${\cal B}^+$ of $T$'s of even cardinality is isomorphic to
$J_2(C)$: every $\eta\in J_2(C)$ can be represented by a divisor $D\sim
kg^1_2-\sum_{i\in T}P_i$, with $|T|=2k$.

\beginlem(3.9)
Let $\eta\in J_2(C)$ be represented by
$D\sim kg^1_2-\sum_{i\in T}P_i$, with $T\in{\cal B}^+$.
Let $S=\P(\pi_*L)$. Then $S_\eta =
\P(\pi_*(L\otimes \eta))$ is isomorphic to $\elm_T(S)$.
\endlem

\bewijs
Factorise $\elm_T$  as  $S\smash{\mathop{\longleftarrow}\limits^{\sigma}}
\wt S  \mapright{\sigma'}S'=\elm_T(S)$: as before, the exceptional
divisors are $E_i$ and $E_i'$. We have $L=E_0+af$; let $T_0\subset T$ be
the index set of points $P_i$ on $E_0$, with cardinality $m$.  Let $\wl
E_0$ be the strict transform of $E_0$ on $\wt S$, and $E_0'$ that on
$S'$.
Then:
$$
\eqalign{
\sigma^*L
&=
\wl E_0+\sum_{j\in T_0}E_j+  af\cr
&=
(\sigma')^*E_0'
-\sum_{i\notin T_0}E_i'+\sum_{j\in T_0}E_j+  af\cr
&\sim
((\sigma')^*E_0'+(a+m-k)f)+\sum_T E_i-kf.\cr}
$$
On $S'$ the linear system  $|E_0'+(a+m-k)f|$ cuts out on $C$ the series
$|L+D |$. We have $E_0'\cdot E_0'=-e+2(k-m)$, so if $2(k-m)< e$, the
section $E_0'$ is the unique negative section, and we see directly that
$E_0'+(a+m-k)f$ is the hyperplane class ($2a=3g+3+e$). Otherwise there is
a section $E_-$ with $E_-\cdot E_-=-e'$, $e'\geq0$, and $E_0'\sim
E_-+(k-m-e/2+e'/2)f$. The hyperplane class is $E_-+(3g+3+e')/2f\sim
E_0'+(a+m-k)f$.
\qed

\roep(3.10) Alternative description of the construction.
We have given the surfaces $S$ as ruled surfaces, birationally embedded
with the linear system $|C|=|2E_0+af|$. Castelnuovo  decribes linear
systems in the plane, of curves of degree $g+e+1$ with one $(g+e-1)$-ple
point and $e-1$ infinitely near double points, if $1\leq e \leq g+1$, or
for $e=0$ of curves of degree $g+3$ with a $(g+1)$-ple point and a double
point in a different point \cite{\rfcas}. If there are at least two
finite double points, a standard Cremona transformation will decrease the
degree, for details see the book \cite{\rfco}.

Now given a curve $C$ in the plane of degree $d$  with $\delta$ double
points (finite or not), and with multiplicity $d-2$ at the origin, and
given a set $T$ of Weierstra\ss\ points of cardinality $k$, we form the
curve $C\cup (\cup_{i\in T}L_i)$ of degree $d+k$, where $L_i$ is the line
joining the origin with the Weierstra\ss\ point $P_i$, and we consider
the linear system  of plane curves of degree $d+k$ with multiplicity
$d+k-2$ at the origin, and $\delta+k$ double points in the double points
of $C$ and in the $P_i$, $i\in T$. As the curves of this system do not
intersect $L_i$ outside $C$,  the rational map $\phi_{|C|}(\P^2) \to
\phi_{|C+\sum L_i|}(\P^2)$ is $\elm_T$.

\roep (3.11) Example {\rm((2.13) }continued\/{\rm )}.
Let $X$ be the cone over the hyperelliptic curve $y^2-1+x^6$, embedded
with $6g^1_2$.  The coordinates on $T^1_X(-1)$ used in (2.13) are
different from the ones in the proof of Proposition (3.6); in principle
one can compute the coordinate transformation, but we will not do this
here. We describe two surfaces with $C$ as hyperplane section, and
identify the induced line in $T^1_X(-1)$.  We remark that the
Weierstra\ss\ points are the sixth roots of unity, and  that the group
$\mu_6$ operates on $\cal B$.

We start with the plane curve $x^2y^2z^2-y^6+x^6$. The curves of degree 6
with multiplicity 4 in $O$, and $A_3$ singularities, tangent to $x=0$ and
$y=0$, are:
$$
z_i=y^{6-i}x^i, \quad \hbox{$i=\mdot06$;}\qquad
w_i=zyxy^{3-i}x^i, \quad \hbox{$i=\mdot03$;}\qquad u=z^2y^2x^2.
$$
This gives equations:
$$
\displaylines{
\Rank \pmatrix{
z_0&\ldots&z_5&w_0&w_1&w_2\cr
z_1&\ldots&z_6&w_1&w_2&w_3\cr} \leq1,\cr
\phi_m=w_{\ab m}w_{\auf m}-uz_{m} , \qquad
\hbox{$m=\mdot0{6}$.}  \cr}
$$
The image of $\P^2$ is isomorphic to the cone over the rational normal
curve of degree 3.  With $t=u-z_0+z_6$ we get a deformation of $C$, with
$s_i=0$.

Now consider the Weierstra\ss\ point $P_0=(1:1:0)$, and the linear
system $|C+L_0|$ of septic curves with additional double point at
$P_0$. A basis is given by:
$$
\displaylines{
\zeta_i=(x-y)^2y^{5-i}x^i, \qquad \hbox{$i=\mdot05$;}
\qquad u_0=z^2y^3x^2, \quad u_1=z^2y^2x^3;\cr
\qquad w_i=zyx(x-y)y^{3-i}x^i, \qquad \hbox{$i=\mdot03$.}
\cr}
$$
On $C$ we have $(y-x)^{-1}\zeta_i= z_{i}-z_{i+1}$, and
$(y-x)^{-1}u_0= z_0+\dots+z_5$, $(y-x)^{-1}u_1=z_1+\dots+z_6$;
therefore we find as equalities on $C$:
$$
\eqalign{
6(y-x)z_0&=u_0+5\zeta_0+4\zeta_1+3\zeta_2+2\zeta_3+\zeta_4\cr
6(y-x)z_1&=u_0-\zeta_0+4\zeta_1+3\zeta_2+2\zeta_3+\zeta_4\cr
6(y-x)z_2&=u_0-\zeta_0-2\zeta_1+3\zeta_2+2\zeta_3+\zeta_4\cr
6(y-x)z_3&=u_0-\zeta_0-2\zeta_1-3\zeta_2+2\zeta_3+\zeta_4\cr
6(y-x)z_4&=u_0-\zeta_0-2\zeta_1-3\zeta_2-4\zeta_3+\zeta_4\cr
6(y-x)z_5&=u_0-\zeta_0-2\zeta_1-3\zeta_2-4\zeta_3-5\zeta_4\cr
6(y-x)z_6&=
u_0-\zeta_0-2\zeta_1-3\zeta_2-4\zeta_3-5\zeta_4-6\zeta_5.\cr
}
$$
Furthermore, $C$ is given by $u_0-u_1-\sum \zeta_i=\tau$. The matrix for
the scroll $\elm_P (S)$ has as first entry $6(y-x)z_0$; i.e.~we make a
coordinate change $\rho z_0:=
u_0+5\zeta_0+4\zeta_1+3\zeta_2+2\zeta_3+\zeta_4$  (we are allowed to
multiply homogeneous coordinates with a common factor).  Below
$6(y-x)z_0$ we have to write
$6(y-x)z_1=u_1+5\zeta_1+4\zeta_2+3\zeta_3+2\zeta_4+\zeta_5$;  the second
entry on the first row is
$u_0-\zeta_0+4\zeta_1+3\zeta_2+2\zeta_3+\zeta_4=6(y-x)z_1+\tau$.
Therefore $s_1=\tau$; proceeding in this way we find $s_i=\tau$ for
$i=\mdot15$. For the coordinates $w_i$ we  take $6(y-x)zyxy^{3-i}x^i$.
Finally we have to compute $t$ in terms of $\tau$. This is now a simple
question of making the explicitly given coordinate transformation in the
equations, and the best way to do this is by computer.  One finds
$t=-4\tau$. Indeed, $s_i=\tau$ and $t=-4\tau$ is a solution of the
equations for the base space.

The other components are found with similar computations.

\roep (3.12)The Milnor fibre.
In their paper \cite{\rflw} Looijenga and Wahl introduce a collection of
smoothing data for singularities, which in many cases distinguish between
components. We determine these data in our case.

Because  the cone over an hyperelliptic curve $C$ of degree $4g+4$, is a
homogeneous singularity, the interior of the Milnor fibre of a one
parameter smoothing, given by  a surface $S$ with $C$ as hyperplane
section, is isomorphic to the open surface $S\setminus C$. We have
described all smoothing components by  smooth ruled surfaces,  except
when $\wl S$ is isomorphic to the cone over a rational normal curve of
degree $g+1$; then $N_{C/S}= (2g+2)g^1_2$. In the last case there does
not exist a one-parameter smoothing, which lies totally in $S^-$, but by
slightly changing the normal bundle we find a smooth ruled surface $S$ on
the same smoothing component, and the Milnor fibre is again $S\setminus
C$.

\beginlem(3.13)
Let $C$ be a curve of type $2E_0+(g+1+e)f$ on $S\cong F_e$. Then the
Milnor fibre $F$ of the cone $X$ over $C$ has homology $H_2(F,{\Bbb
Z})=\Z^{2g+1}$ and  $H_1(F)=\Z/2$. If $g+e+1\cong 0 \pmod 2$, then
$H_1(F)=\Z/2$ and $H^2(F,\partial F)=\Z^{2g+1}\oplus \Z/2$; otherwise
$H_1(F)=0$ and $H^2(F,\partial F)=\Z^{2g+1}$.
\endlem

\bewijs
We compute the homology of $F=S\setminus C$ with the exact homology
sequence of the pair $(S,S\setminus C)$. For the surface $S$ we have
$H^3(S)=H^1(S)=0$ and $H_2(S)=\Pic (S)$. We may  take a new basis,
consisting of $E=E_0+\ab{g+1+e}f$ and $f$. Then $C\sim 2E+\epsilon f$
with $\epsilon=0$ if $g+1+e$ is even, and $1$ otherwise, and $E\cdot
E=g+1-\epsilon$.   We have already encountered this division in cases.
Let $N$ be a tubular neighbourhood of $C$;  by excision and Poincar\'e
Duality $H_i(S,S\setminus C)=H_i(N,\partial N)\cong H^{4-i}(C)$. Therefore
we have:
$$
\cdmatrix{
 H_3(S,S\setminus C)& \hookrightarrow & H_2(S\setminus C)&
 \to& H_2(S)& \to &
 H_2(S,S\setminus C) &\twoheadrightarrow& H_1(S\setminus C) \cr
\Vert\wr \atop\displaystyle {\Bbb Z}^{2g}&&&&
 \Vert\wr \atop\displaystyle {\Bbb Z}^2&
 \searrow&\mapdown {PD}&&\cr
&&&&&&H^2(C) && \cr}
$$
The composed map $H_2(S) \to H^2(C)$ is given by intersecting with the
curve $C$; so $f\mapsto 2$ and $E\mapsto 2g+2-\epsilon$. If $\epsilon=0$,
then $H_1(F;{\Bbb Z})={\Bbb Z}/2$, and  $H_1(F;{\Bbb Z})=0$, if
$\epsilon=1$. In the odd case  $\ker \{H_2(S)\to H_2(S,S\setminus C)\}$
is generated by $2E-(2g+2-\epsilon)f$, and
$(2E-(2g+2-\epsilon)f)^2=-4(g+1)$; in the even case a generator is
$E-(g+1)f$, with self-intersection $-(g+1)$.

We compute $H_2(F,\partial F)\cong H^2(F)$ with the cohomology sequence
of the pair $(S,S\setminus C)$. We obtain:
$$
0\longrightarrow H^2(S,S\setminus C)
 \longrightarrow H^2(S) \longrightarrow
 H^2(S\setminus C) \longrightarrow H^3(S,S\setminus C)
 \longrightarrow 0.
$$
By Poincare duality $H^i(S,S\setminus C;{\Bbb Z})\cong H_{4-i}(C)$, and
the map $H_2(C)\to H^2(S)$ is the transpose of the map $H_2(S) \to
H^2(C)$ above.
\qed

The boundary $\partial F$ of the Milnor fibre is diffeomorphic to the
link $M$ of the singularity. As $M$ is a circle bundle over $C$, we have
that $H_1(M)=\Z^{2g}\oplus\Z/(4g+4)$. The linking form on $H_1(M)_t$, the
torsion part, lifts to a quadratic function, which can be computed from
the resolution \cite{\rflw, 4.2, 4.6}. We have, with $\bar e$ a generator
of $H_1(M)_t$:
$$
q\colon H_1(M)_t \to {\Bbb Q}/\Z, \qquad  q(m\bar e) = {m(m-6g-2)\over
8g+8}.
$$
By \cite{loc.~cit., Thm.~4.5} the group $H_2(F)$, with intersection form
$\cal I$, carries the structure of a quadratic lattice with associated
non-degenerate lattice $\wl H_2(F)=H_2(F)/{\rm rad}\cal I$, whose
discriminant quadratic function is canonically isomorphic to
$(I^\perp/I,q_I)$, where $I$ is the $q$-isotropic subgroup ${\rm Im} \{
\partial_t\colon H_1(F,M)_t\to H_1(M)_t\}$.

In our case the radical of the intersection form has rank $2g$:
$\mu_0=2g$, $\mu_-=1$ and $\mu_+=0$.  If $g+1+e$ is odd, then
$H_1(F,M)_t=0$, so $I=0$ and $I^\perp =H_1(M)_t$. For  $g+1+e$ even $\wl
H_2(F)$ is generated by $E-(g+1)f$, with self-intersection $-(g+1)$. In
this case the isotropic subgroup $I$ has order two.

We recall the definition of smoothing data \cite{loc.~cit., 4.16} in the
case $\mu_+=0$. Let $A$ be a finitely generated abelian group, and
$q\colon A_t\to {\Bbb Q}/\Z$ a nonsingular quadratic function on its
torsion part. Then ${\cal S}(A,q)$ is the set of equivalence classes of
4-tuples $(V,Q,I,i)$, where $(V,Q)$ is a negative semi definite ordinary
lattice, $I$ is a $q$-isotropic subspace of $A_t$, and  $i\colon
V^*/B'(V)\to A/I$ is an injective homomorphism with finite cokernel,
which induces an isomorphism $\wl V^\#/\wl V\to i^\perp/I$ (here
$B'\colon V\to V^*$ is the adjoint of the bilinear part of $Q$); two
4-tuples $(V_1,Q_1,I_1,i_1)$ and $(V_2,Q_2,I_2,i_2)$ are equivalent,  if
$I_1=I_2$ and there exists an isomorphism
$\Phi\colon(V_1,Q_1)\to(V_2,Q_2)$ such that $i_2=\phi\circ i_1$, where
$\phi$ is induced by $\Phi$.

Let $J$ be the preimage of ${\rm Im }(i)$ in $A$. Then $J\cap
A_t=I^\perp$ (so $J$ determines $I$). By \cite {loc.~cit., 4.17} the
elements in ${\cal S}(A,q)$, which give the same triple $(V,Q,J)$ form
$O(q_I)$-orbit, where $O(q_I)$ is the group of linear automorphims
preserving the quadratic function $q_I$ on $I^\perp/I$. Each orbit is
equivalent to a coset space in $O(q_I)$ of  ${\rm Im} \{O(\wl Q)\to
O(q_I)\}$.

For all our smoothings we have $H_1(M)_t/I^\perp \cong H_1(F)$. In terms
of $J$ this means that $A_t+J=A$. In \cite{loc.~cit.} this condition is
shown to hold for all minimally elliptic singularities, with the concept
of permissible quotients. Denote for our singularity $X$ by ${\cal S}(X)$
the subset of ${\cal S}(H_1(M),q)$ of lattices with $\mu_-=1$,
$\mu_0=2g$, for which $H_1(M)_t+J=H_1(M)$.  The set of smoothing
components maps to ${\cal S}(X)$.

\beginprop(3.14)
Suppose $g$ is even. Then the set ${\cal S}(X)$ has $1+2^{2g}$ elements.
On the set of even smoothing components {\rm (}i.e. with $g+1+e$
even\/{\rm )} the map to ${\cal S}(X)$ is injective, whereas the odd
components all map to the same element.
\endprop

\bewijs
We first do not impose  restrictions on $g$. By \cite{\rflw, Ex. 4.6} the
group $H_1(M)_t$ has an isotropic subgroup of order $r$, if $4g+4=r^2s$,
$2g-2=ru$ with $s(1+r)\equiv u \pmod 2$. So $8\equiv 0 \pmod r$. Possible
values are $r=2$, $u=g-1$, $s=g+1$, for arbitrary $g$. If $g=8k+3$, then
also $(r,s,u)=(4,2k+1,4k+1)$ is a solution.

Next we determine $O(q_I)$. Write $g+1=h$. First consider the case $I=0$.
We look for a map $\bar e\mapsto d \bar e$ such that for all $m$:
$$
8h(q(md\bar e)-q(m\bar e)=(d-1)(m^2(d+1)-6hm+4m) \equiv 0 \pmod{8h}.
$$
Let $G=\gcd(d-1,8h)$, and $8h=aG$, then  $m^2(d+1)-6hm+4m\equiv 0 \pmod
a$, and therefore $d+1\equiv 6h-4$, $2(d+1)\equiv 0$, so $16\equiv 0
\pmod a$. Because  $d<4h$, possible values for $a$ are 4, 8 and 16. Only
for $a=8$ we find  solutions, given by $d=h+1$, $h=8k-2$ or $d=3h+1$,
$h=8k+2$. In these cases $O(q_I)$ is a group of order two. If $I$ is a
group of order 2, with $2h\bar e$ as generator, then $I^\perp/I\cong
\Z/h$, and we find that for $h=8k$, $d=4k+1$, the group $O(q_I)$ has
order two.

So if $h$ is odd, then the group $O(q_I)$ is trivial, and the number of
elements ${\cal S}(X)$ is the number of triples $(V,Q,J)$. If $I=0$, then
$J=H_1(M)_t$, and this is the only possibility. Now suppose that $I$ is
a group of order two. In \cite{loc.~cit., 6.2} it is shown that $\Hom
(H_1(M)/ H_1(M)_t,H_1(M)_t/I^\perp)$ acts simply transitive on the
collection of subgroups $J$ of $H_1(M)$ with $H_1(M)_t\cap J=I^\perp$ and
$H_1(M)_t+J=H_1(M)$. In our case this group is isomorphic to
$(\Z/2)^{2g}$.

To show injectivity, we have to identify this action with the action of
the ${\cal B}^+$.  The group $H_2(F,\partial F)$ is isomorphic to
$H_2(S,C)$; we descibe this group with the exact sequence of the pair
$(S,C)$. For $H_2(S)$ we have two generators, $E$ and $f$, and $C\sim
2E$. Take a basis of $H_1(C)$, in the traditional way (cf. \cite{\rfTata,
p. 3.75}): choose a certain system of $2g$  paths in $\P^1-B$, $B$ the
set of branch points of $\pi\colon C\to\P^1$, which each encircle an even
number of branch points, and lift them to $C$; in $\P^1$ such a path
$\alpha_i$ bounds a disc $D_i$. Now we consider the surface
$S\mapright{\pi} \P^1$. In $S$ a lift of $\alpha_i$ to $C$ is the
boundary of a lift $\Delta_i$ of $D_i$ to $S$. Then $E$, $f$ and the
$\Delta_i$ generate $H_2(F, \partial F)$. The map to $H_1(M)$ is given by
intersecting on $S$ with the boundary $M$ of a tubular neighbourhood of
$C$ (with as orientation the one induced by $F$). The manifold $M$ is a
circle bundle over $C$, with Euler number $-(4g+4)$; let $\gamma$ be the
homology class of a fibre, then $\gamma$ generates  $H_1(M)_t$. A general
fibre of $S\mapright{\pi} \P^1$ intersects $M$ in two circles; let
furthermore $A_i$ be a lift of $\alpha_i$ to $M$. Then the map $\partial
\colon H_2(F,\partial F)\to H_1(M)$ is given by $\partial(f)= 2\gamma$,
$\partial (E)=2(g+1)\gamma$ and $\partial (\Delta_i)=A_i$.

Consider now $S'=\elm_T(S)$, giving the Milnor fibre $F'$. Let $\alpha_i$
encircle a point $P\in T$; we may suppose that $\Delta_i$ does not pass
through $P$, but intersects the fibre $\pi^{-1}(P)$ once, outside the
tubular neighbourhood of $C$..  The strict transform $\Delta_i'$ on $S'$
is the element of $H_2(F',\partial F')$, which lifts $\alpha_i$. But on
$S'$ the disc $\Delta_i$ intersects $C'$ transversally, and therefore
intersects $M$ in a circle. So $\partial ' (\Delta_i')=A_i+t_i\gamma$,
where $t_i$ is the number of branch points lying in $D_i$. If $t_i$ is
odd,  then  $\partial ' (\Delta_i')$ does not lie in  $\partial (H_2F,
\partial F)$. This construction establishes an isomorphism between ${\cal
B}^+$ and $\Hom(H_1(M)/H_1(M)_t, H_1(M)_t/I^\perp)$.
\qed

\roep (3.15)Remark.
In the case $g=3$ an isotropic subgroup $I\subset H_1(M)_t$ of order
4 occurs for a smoothing of $X(C,K^4)$: a non hyperelliptic curve is
hyperplane section of $4$-fold Veronese embedding of $\P^2$, i.e.
$\P^2$ embedded with $\sier {}(4)$. It is well known that
$H_1(\P^2\setminus C)=\pi_1(\P^2\setminus C)=\Z/4$.

\roep (3.16)Remark. For a simple elliptic singularity of multiplicity
$8=4g+4$ the topological computation still work; in this case for $I=0$
the group $O(q_I)$ is trivial. Then the permissible smoothing data form a
set of $1+2^2=5$ elements \cite{\rflw, 6.4}. The number of smoothing
components is also five; there are four `even' components and one `odd'
component.

Let $X$ be the cone over an elliptic curve $E$ of degree 8, embedded with
the linear system $L$; so $\phi_L\colon E\to\P^7$.  We identify
$Pic^8(E)$ with ${\rm Jac}(E)\cong E$ in such a way that $L$ represents
the origin of the group law: we write $L=8P$. There are 4 line bundles
$M$ of degree 4 with $M^{\otimes 2}=L$. Each of these embeds $E$ in
$\P^3$; then $\phi_L$ is the composition of $\phi_M$ with the  the
Veronese embedding  $V\colon \P^3\to\P^9$. The curve $\phi_M(E)$ lies on
a pencil of quadrics, which is transformed by $V$ in a pencil of
hyperplane sections of $V(\P^3)$; the common intersection is $E$. This
construction gives  a 2-dimensional linear subspace (a cone over $\P^1$)
in the negative part of the base space of the versal deformation of $X$.
The general quadric in the pencil is non-singular and has two rulings,
but there are four singular quadrics, which are cones with one ruling; the
resulting double cover of $\P^1$ is isomorphic to $E$. The fifth
component is isomorphic to the cone over $E$;  given a point $Q\in E$, we
embed  $E$ in the plane with the linear system $|2P+Q|$, then  embed
$\P^2$ with the triple Veronese embedding $V_3$ in $\P^9$, and project
the surface $V_3(\P^2)$ from the point $V_3\circ\phi_{|2P+Q|}(3Q)$ to
$\P^8$. Each of the resulting surfaces has a unique ruling.

Each ruling determines a $g^1_2\in \Pic^2\cong E$; on the even components
a given $g^1_2$ occurs exactly once, and on the odd component four times.
The operations $\elm_T$, with $T$ a subset of the set of branch points
$\cal B$ of the $g^1_2$, permute the $2^3$ points in the projectivised
base space.

\roep(3.17) Other degrees.
The structure of the equations (see (2.2) and (2.12)) and the dimension
of $\ttw{-2}$, which is $d-2g-3$, two less than the number of equations,
which cut out the curve on the scroll, lead us to expect that the base
space $S^-$ in negative degree is a complete intersection of dimension
$4g+5-d$, if $2g+4\leq d \leq 4g+5$. The number of equations will
probably increase linearly for $d>4g+5$, until it is $(g+1)(2g+3)$, the
number of quadratic monomials in $2g+2$ variables.

For $d=4g+5$ the cone $X$ has only conical deformations \cite{\rftend},
and therefore all infinitesimal deformations of negative degree are
obstructed. In this case we have indeed a zero-dimensional complete
intersection.

For $d=4g+3$ the space $S^-$ should be the cone over a complete
intersection curve of degree $2^{2g}$, which is a $2^{2g}$-fold
unbranched covering of $C$. One finds this degree from the genus
formula's. From the point of view of the surfaces, which correspond to
one parameter deformations, we get the following picture. A surface $S$
of degree $4g+3$ is a scroll, blown up in one point $P$ of the
hyperelliptic curve $C$. The fibration $\pi\colon S\to\P^1$ has one
singular fibre, which intersects $C$ in $P$ and in $\wl P$, the image of
$P$ under the hyperelliptic involution. By blowing down one of the two
$(-1)$-curves in the exceptional fibre we get a ruled surface with a
minimal section $E_0$ with $E_0^2=-e$; if we blow down the other curve,
then the parity of $e$ is just the opposite. Suppose that blowing down the
curve through $P$ gives an even $e$. As we have seen, there are $2^{2g}$
different surfaces with $e$ even, such that the characteristic linear
system on $C$ is a given linear system $|L+P|$. Altogether we get a
covering of $\P^1$ of degree $2^{2g+1}$, with simple branching at  the
Weierstra\ss\ points, so indeed  a $2^{2g}$-fold unbranched covering of
$C$. In particular, there is only one smoothing component, except of
course when $g=3$.

We leave it to the reader to continue this description for still lower
values of the degree $d$.


\vskip0pt plus.3\vsize\penalty-250
 \vskip0pt plus-.3\vsize
\beginsection Bibliography

\frenchspacing
\parindent=45pt
\bigskip
\def\item#1{\par\noindent\rlap{#1}\hskip\parindent%
\hangindent\parindent\ignorespaces}

\item{\EGA}
  A. Grothendieck,
  {\sl \'El\'ements de g\'eom\'etrie alg\'ebrique.\/} R\'edig\'es avec
  la collaboration de J. Dieudonn\'e. IV, {\sl \'Etude locale des
  sch\'emas et des morphismes de sch\'emas\/.} Publ. Math. IHES
  {\bf 20} (1964), {\bf 24} (1965), {\bf 28} (1966), {\bf 32} (1967).

\medskip
\parindent=0pt
\def\item{\par\hangindent=20pt}

\item
  Enrico Arbarello e Maurizio Cornalba,
  {\sl Su una congettura di Petri.\/}
  Comment. Math. Helv. {\bf 56}, (1981), 1--38.

\item
  Enrico Arbarello and Eduardo Sernesi,
  {\sl Petri's Approach to the Study of the Ideal Associated to a
  Special Divisor.\/}
  Invent. Math. {\bf 49}, (1978), 99--119.

\item
	Dave Bayer and Mike Stillman,
	{\sl Macaulay: A system for computation in
	algebraic geometry and commutative algebra.\/}
	Source and object code available for Unix and Macintosh
	computers. Contact the authors, or download from
	{\tt  zariski.harvard.edu} via anonymous ftp.

\item
  Kurt Behnke and Jan Arthur Christophersen,
  {\sl Hypersurface sections and obstructions (rational surface
  singularities).\/}
  Comp. Math. {\bf 77} (1991), 233--268.

\item
  Guido Castelnuovo,
  {\sl Sulle superficie algebriche le cui sezioni piane sono curve
  iperellitiche\/.}
  Rend. Circ. Mat. Palermo {\bf 4} (1890), 73--88.

\item
  Ciro Ciliberto, Joe Harris and Rick Miranda
  {\sl On the surjectivity of the Wahl map.\/}
  Duke J. Math. {\bf 57} (1988), 829--858.

\item
  Fabio Conforto,
  {\sl Le superficie razionale.\/}
  Bologna, Zanichelli, 1939.

\item
   R. Drewes,
   {\sl Infinitesimale Deformationen von Kegeln \"uber
   trans\-ka\-no\-nisch eingebetteten hyperelliptischen Kurven.\/}
   Abh. Math. Sem. Univ. Hamburg {\bf 59} (1989), 269--280.

\item
  Mark L. Green,
  {\sl Koszul cohomology and the geometry of projective varieties.\/}
  J. Diff. Geom. {\bf 19} (1984), 125--171.

\item
  Gert-Martin Greuel,
  {\sl On deformations of curves and a formula of De\-ligne.\/}
  In: Algebraic Geometry, La R\'abida 1981, pp. 141--168.
  Berlin etc., Springer 1982 (Lect. Notes in Math.; 961).

\item
  Robin Hartshorne,
  {\sl Algebraic Geometry.\/}
  Berlin etc., Springer 1977.

\item
  Steven L. Kleiman,
  {\sl The Enumerative Theory of Singularities.\/}
  In:  Real and Complex Singularities, Oslo 1976, pp. 297--396
  Alphen a/d Rijn, Sijthoff \&~ Noordhoff 1977.

\item
  Paulette L\'evy-Bruhl-Mathieu.
  {\sl \'Etude des surfaces d'ordre\/ $p - 2 + h$
  $(0\leq h<p)$ passant par une courbe canonique de genre
  $p$. Application \`a la classification des courbes
  al\'gebriques de genre inf\'erieur \`a 11.}
  C. R. Ac. Sc. Paris, 236 (1953), 2032- 2034.

\item
  Eduard Looijenga,
  {\sl The smoothing components of a triangle singularity II.\/}
  Math. Ann. {\bf 269} (1984), 357--387.

\item
  Eduard Looijenga and Jonathan Wahl,
  {\sl Quadratic functions and  smoothing surface singularities.\/}
  Topology {\bf 25} (1986), 261--291.

\item
  David Mumford,
  {\sl A remark on the paper of M.~Schlessinger.\/}
  Rice Univ. Studies {\bf 59} (1973), 113-117.

\item
  David Mumford,
  {\sl Tata Lectures on Theta II.\/}
  Boston etc., Birkh\"auser 1984 (Progress in Math.; 28).

\item
  Ragni Piene,
  {\sl Polar classes of singular varieties.\/}
  Ann. Sc. \'Ec. Norm. Sup.  {\bf 11} (1978) 247--276.

\item
  Henry C. Pinkham,
  {\sl Deformations of algebraic varieties with $G_{m}$-action.\/}
  Ast\'erisque {\bf 20} (1970).

\item
  Miles Reid,
  {\sl Surfaces with $p_g=3, K^2=4$ according to E. Horikawa
  and D. Dicks.\/}
  Text of a lecture, Univ. of Utah  and Univ. of Tokyo 1989.

\item
  Michael Schlessinger,
  {\sl On rigid singularities.\/}
   Rice Univ. Studies {\bf 59} (1973), 147-162.

\item
  Frank-Olaf Schreyer,
  {\sl Syzygies of Canonical Curves and Special Linear Series.\/}
  Math. Ann. {\bf 275} (1986), 105--137.

\item
  Francesco Severi,
  {\sl Vorlesungen \"uber Algebraische Geometrie.\/}
  Leipzig, Berlin, B. G. Teubner, 1921.

\item
  Jan Stevens,
  {\sl On the number of points determining a canonical curve.\/}
  Indag. Math. {\bf  51} (1989), 485--494.

\item
  Jan Stevens,
  {\sl On deformations of singularities of low dimension and
  codimension: base spaces and smoothing components.\/}
  Preprint 1991 (Schriften\-reihe Forschungs\-schwerpunkt Komplexe
  Mannigfaltig\-kei\-ten; 129).

\item
  Sonny Tendian,
  {\sl Deformations of Cones and the Gaussian-Wahl Map.\/}
  Preprint 1990.

\item
  Sonny Tendian,
  {\sl Surfaces of degree $d$ with sectional genus $g$ in
  $\P^{d+1-g}$ and deformations of cones.\/}
  Duke J. Math. {\bf 65} (1992), 157--185.

\item
  Sonny Tendian,
  {\sl Extensions of hyperelliptic curves.\/}
  Preprint 1992.

\item
  Jonathan M. Wahl,
  {\sl The Jacobian algebra of a graded Gorenstein Singularity.\/}
  Duke J. Math. {\bf 55} (1987) 843--871.

\item
  Jonathan Wahl,
  {\sl Deformations of Quasi-homogeneous Surface Singularities.\/}
  Math. Ann. {\bf 280} (1988), 105--128.

\item
  Jonathan Wahl,
  {\sl Introduction to Gaussian maps on an algebraic curve.\/}
  Notes prepared in connection with lectures at the Trieste
  Conference on Projective Varieties, June 1989.

\item
  Jonathan Wahl,
  {\sl Gaussian maps on algebraic curves.\/}
  J. Diff Geom. {\bf 32} (1990), 77--98.


\vfill
\parindent=0pt
{\obeylines
Address of the author:
Mathematisches Seminar der Universit\"{a}t Hamburg
Bundesstra\ss e 55, D 2000 Hamburg 13, Germany}
\settabs\+E-mail: &(new) &\cr
\+E-mail: &(new) & stevens@geomat.math.uni-hamburg.de\cr
\+& (old) & ms00010@dhhuni4.bitnet\cr

\bye